\documentclass[11pt]{article}
\usepackage{mathtools,amsmath,amsthm,amssymb,mathrsfs,amsfonts}
\usepackage{cases}
\usepackage{leftidx}
\usepackage{epsfig}
\usepackage{leftidx}
\usepackage{graphicx,subfigure}
\usepackage{color, epstopdf}
\usepackage{cite}
\usepackage{graphicx,wrapfig,tabularx}
\usepackage{flafter}
\usepackage{fancyhdr}
\usepackage{stmaryrd}
\usepackage{multicol,multirow,booktabs}
\usepackage{dsfont}
\usepackage{booktabs,threeparttable}
\usepackage[center]{caption2}
\usepackage{epstopdf}
\usepackage [latin1]{inputenc}
\usepackage{multirow}
\usepackage{enumerate}
\usepackage{enumitem}
\usepackage{algpseudocode,algorithm,algorithmicx}

\newtheorem{theorem}{Theorem}[section]
\newtheorem{lemma}{Lemma}[section]
\newtheorem{proposition}{Proposition}[section]

\newtheorem{example}{Example}
\newtheorem{remark}{Remark}

\def\ck{{\scriptstyle{K}}}
\def\hck{{\hat{\scriptstyle{K}}}}
\def\Ck{{{\mathtt{K}}}}

\newcommand{\myvec}[1]{\boldsymbol{#1}}

\newcommand{\zd}{\,\mathrm{d}}

\newcommand{\abs}[1]{\left|#1\right|}
\newcommand{\absb}[1]{\big|#1\big|}

\newcommand{\abst}[1]{|#1|}
\newcommand{\bra}[1]{\left(#1\right)}
\newcommand{\brab}[1]{\big(#1\big)}
\newcommand{\braB}[1]{\Big(#1\Big)}
\newcommand{\brat}[1]{(#1)}
\newcommand{\kbra}[1]{\left[#1\right]}
\newcommand{\kbrab}[1]{\big[#1\big]}
\newcommand{\kbraB}[1]{\Big[#1\Big]}

\newcommand{\myinner}[1]{\left\langle#1\right\rangle}
\newcommand{\myinnert}[1]{\langle#1\rangle}
\newcommand{\myinnerb}[1]{\big\langle#1\big\rangle}
\newcommand{\myinnerB}[1]{\Big\langle#1\Big\rangle}
\newcommand{\mynorm}[1]{\left\|#1\right\|}
\newcommand{\mynormb}[1]{\big\|#1\big\|}

\newcommand{\mynormt}[1]{\|#1\|}

\title{A class of refined implicit-explicit Runge-Kutta methods\\ with robust time adaptability and unconditional convergence\\ for the Cahn-Hilliard model}
\author{Hong-lin Liao\thanks{ORCID 0000-0003-0777-6832. School of Mathematics,
	Nanjing University of Aeronautics and Astronautics,
	Nanjing 211106, China; Key Laboratory of Mathematical Modeling and High Performance Computing of Air Vehicles (NUAA), MIIT, Nanjing 211106, China. Emails: liaohl@nuaa.edu.cn and liaohl@csrc.ac.cn.
	This author is supported by NSF of China under grant numbers 12471383 and 12071216.}
	\quad Tao Tang\thanks{School of Electrical and Computer Engineering, Guangzhou Nanfang College, and Institute for Advanced Study, BNU-HKBU United International College, China. Email: {ttang@uic.edu.cn}. This author is supported by NSF of China under grant numbers 11731006 and K20911001.}
	\quad Xuping Wang\thanks{School of Mathematics, Nanjing University of Aeronautics and Astronautics, Nanjing 211106, China. Email: wangxp@nuaa.edu.cn.}
	\quad Tao Zhou\thanks{Institute of Computational Mathematics and Scientific/Engineering Computing,	Academy of Mathematics and Systems Science, Chinese Academy of Sciences, Beijing, 100190, China. Email: {tzhou@lsec.cc.ac.cn}. This author is supported by NSF of China under grant number 12288201.}
}

\begin{document}
	
	\maketitle
	
	\begin{abstract}
		One of main obstacles in verifying the energy dissipation laws of implicit-explicit Runge-Kutta (IERK) methods for  phase field equations is to establish the uniform boundedness of stage solutions without the global Lipschitz continuity assumption of nonlinear bulk. With the help of discrete orthogonal convolution kernels, an updated time-space splitting technique is developed to establish the uniform boundedness of stage solutions for a refined class of IERK methods in which the associated differentiation matrices and the average dissipation rates are always independent of the time-space discretization meshes. This makes the refined IERK methods highly advantageous in self-adaptive time-stepping procedures as some larger adaptive step-sizes in actual simulations become possible. From the perspective of optimizing the average dissipation rate, we construct some parameterized refined IERK methods up to third-order accuracy, in which the involved diagonally implicit Runge-Kutta methods for the implicit part have an explicit first stage and allow a stage-order of two such that they are not necessarily algebraically stable. Then we are able to establish, for the first time, the original energy dissipation law and the unconditional $L^2$ norm convergence. Extensive numerical tests are presented to support our theory.
		\\[1ex]		
			\\[1ex]
	\end{abstract}

\section{Introduction}
\label{sec: introduction}
\setcounter{equation}{0}

The aim of this work is to present a class of refined implicit-explicit Runge-Kutta methods with robust time adaptability and unconditional convergence for the Cahn-Hilliard model. In our previous work \cite{LiaoWangWen:2024arxiv}, a unified theoretical framework was suggested to examine the energy dissipation properties at all stages of additive implicit-explicit Runge-Kutta (IERK) methods up to fourth-order accuracy for gradient flow problems. Some parameterized IERK methods preserving the original energy dissipation law unconditionally were constructed by applying the so-called first same as last method, that is, the diagonally implicit Runge-Kutta (DIRK) method with the explicit first stage and stiffly-accurate assumption for the linear stiff term, and applying the explicit Runge-Kutta method for the nonlinear term. The main idea in \cite{LiaoWangWen:2024arxiv} is to construct the differential forms and the associated differentiation matrices of IERK methods by using the difference coefficients of method and the so-called discrete orthogonal convolution kernels \cite{LiaoJiWangZhang:2022,LiaoJiZhang:2022PFC,LiaoZhang:2021}. Nonetheless, this theory was built on a strong assumption: the involved nonlinear terms are globally Lipschitz continuous.

In this work, we will consider a much more practical assumption: the nonlinear term is continuously differentiable. Under this assumption, we shall perform a unified $L^2$ norm convergence analysis of a refined class of IERK methods and establish the associated original energy dissipation laws for the well-known Cahn-Hilliard (CH) equation. To the best of our knowledge, this is the first time such original energy dissipation law and unconditional $L^2$ norm convergence of IERK methods are established for the CH model without the global Lipschitz continuity assumption of the nonlinear bulk.

Consider a free energy functional of Ginzburg-Landau type,
\begin{align}\label{cont:free energy}
	E[\Phi] = \int_{\Omega}\kbra{\tfrac{\epsilon^2}{2}\abst{\nabla\Phi}^2+F(\Phi)}\zd\myvec{x}\quad\text{with}\quad F(\Phi):=\tfrac14(\Phi^2-1)^2,
\end{align}
where $\myvec{x}\in\Omega\subseteq\mathbb{R}^2$ and $0<\epsilon<1$ is proportional to the interface width.
The well known Cahn-Hilliard equation is given by the $H^{-1}$ gradient flow
associated with the free energy functional $E[\Phi]$,
\begin{align}\label{cont: Problem-CH}
	\partial_t \Phi=\Delta\kbra{F'(\Phi)-\epsilon^2\Delta\Phi}\quad\text{for $\myvec{x}\in\Omega$}.
\end{align}
Assume that $\Phi$ is periodic over the domain $\Omega$.
By applying the integration by parts, one can find the volume conservation, $\brab{\Phi(t),1}=\brab{\Phi(t_0),1}$,
and the following original energy dissipation law,
\begin{align}\label{cont:energy dissipation}
	\tfrac{\zd{E}}{\zd{t}}=\brab{\tfrac{\delta E}{\delta \Phi},\partial_t \Phi}_{L^2}
	=-\bra{(-\Delta)^{-1}\partial_t \Phi,\partial_t \Phi}_{L^2}\le 0,
\end{align}
where the $L^{2}$ inner product $\bra{u,v}_{L^2}:=\int_{\Omega}uv\zd{\myvec{x}}$ for all $u,v\in{L}^{2}(\Omega)$.

Always, the explicit approximation of the nonlinear bulk $f$ will be adopted in our IERK methods so that they are computationally efficient by avoiding the inner iteration at each time stage. To control the possible instability stemmed from the explicit approximation of $f$, 
we introduce the following stabilized operators with the stabilized parameter $\kappa\ge0$, cf. \cite{DuJuLiQiao:2021SIREV,FuTangYang:2024,LiaoWang:2024MCOM,LiaoWangWen:2024JCP}, 
\begin{align}\label{def: stabilized parameter}
	L_{\kappa}\Phi:=-\epsilon^2\Delta\Phi+\kappa \Phi \quad\text{and}\quad f_{\kappa}(\Phi):=\kappa \Phi-F'(\Phi),
\end{align} 
such that the problem \eqref{cont: Problem-CH} becomes the stabilized version
\begin{align}\label{cont: Problem-CH stabilized}
	\partial_t \Phi=\Delta\kbra{L_{\kappa}\Phi-f_{\kappa}(\Phi)}\quad\text{for $\myvec{x}\in\Omega$}.
\end{align}
Interested readers are referred to the series \cite{LiQiao:CMS2017,LiQiao:JSC2017,LiQiaoTang:SINUM2018,LiQuanTang:MC2022}, where Li et al. characterized the sizes of stabilization parameter for a wide class of semi-implicit stabilized methods for the CH model. 

We apply the Fourier pseudo-spectral method to approximate the spatial operators $\Delta$ and $L_{\kappa}$, as described in Section \ref{sec: motivations}, with the associated discrete operators (matrices) $\Delta_h$ and $L_{\kappa,h}$, respectively. For a finite $T$, consider a nonuniform mesh $0=t_0<t_1<\cdots<t_N=T$ with the time-step $\tau_n=t_n-t_{n-1}$. The time operator is approximated by IERK methods and let $\phi_h^k$ be the numerical approximation of $\Phi_h^k:=\Phi(\myvec{x}_h,t_k)$ at the discrete time level $t_k$ for $0\le k\le N$. For a $s$-stage Runge-Kutta method, let $u_h^{n,i}$ be the approximation of $\Phi_h^{n,i}:=\Phi(\myvec{x}_h,t_{n-1}+c_{i}\tau_n)$ at the abscissas $c_1:=0$, $c_i>0$ for $2\le i\le s-1$, and $c_{s}:=1$.   
To integrate the nonlinear model \eqref{cont: Problem-CH stabilized} from $t_{n-1}$ ($n\ge1$) to the next grid point $t_n$, we consider the following $s$-stage IERK method
\cite{AscherRuuthSpiteri:1997,FuTangYang:2024,LiaoWangWen:2024arxiv,ShinLeeLee:2017,ShinLeeLee:2017CMA}
	\begin{align}	\label{Scheme: general IERK}		
		u_h^{n,i}:=&\,u_h^{n,1}+\tau_n\sum_{j=1}^{i}a_{i,j}\Delta_h L_{\kappa,h}u_h^{n,j}
		-\tau_n\sum_{j=1}^{i-1}\hat{a}_{i,j}\Delta_h f_{\kappa}(u_h^{n,j})\quad\text{for $n\ge1$ and $1\le i\le s$},		
	\end{align}
where $u_h^{n,1}:=\phi_h^{n-1}$ and $\phi_h^{n}:=u_h^{n,s}$. The associated Butcher tableaux reads
\begin{equation*}
	\begin{array}{c|c}
		\mathbf{c} & A \\
		\hline\\[-8pt] 	& \mathbf{b}^T
	\end{array}
	=	\begin{array}{c|ccccc}
		c_{1} & 0 &  &  &  &   \\
		c_{2} & a_{21} & a_{22} &  &  &   \\
		c_{3} & a_{31} & a_{32} & a_{33} &  &   \\
		\vdots & \vdots & \vdots & \ddots & \ddots &  \\[2pt]
		c_{s} & a_{s,1} & a_{s,2} &  \cdots  & a_{s,s-1}   & a_{s,s} \\[2pt]
		\hline  & a_{s,1} & a_{s,2} &  \cdots  & a_{s,s-1}   & a_{s,s}
	\end{array},\quad
	\begin{array}{c|c}
		\mathbf{c} & \widehat{A} \\
		\hline\\[-8pt] 	& \hat{\mathbf{b}}^T
	\end{array}
	=\begin{array}{c|ccccc}
		c_{1} & 0 &  &  &  &   \\
		c_{2} & \hat{a}_{21} & 0 &  &  &   \\
		c_{3} & \hat{a}_{31} & \hat{a}_{32} & 0 &  &   \\
		\vdots & \vdots & \vdots & \ddots & \ddots &  \\[2pt]
		c_{s} & \hat{a}_{s,1} & \hat{a}_{s,2} &  \cdots  & \hat{a}_{s,s-1}   & 0 \\[1pt]
		\hline\\[-8pt]  & \hat{a}_{s,1} & \hat{a}_{s,2} &  \cdots  & \hat{a}_{s,s-1}   & 0
	\end{array}.
\end{equation*}
Without losing the generality, we always assume that $\hat{a}_{k+1,k}(z)\neq0$ for any $1\le k\le s-1$. 
That is, the stiff linear term $L_{\kappa}$ is approximated by the $s$-stage stiffly-accurate DIRK method with the coefficient matrix $A$, the abscissa vector $\mathbf{c}=A\mathbf{1}$ and the vector of weights $\mathbf{b}:=A^T\mathbf{e}_s$, while the nonlinear term $f_{\kappa}$ is approximated by the $s$-stage  explicit Runge-Kutta methods with the strictly lower triangular coefficient matrix $\widehat{A}$, the abscissa vector $\hat{\mathbf{c}}=\widehat{A}\mathbf{1}$ and the vector of weights $\hat{\mathbf{b}}:=\widehat{A}^T\mathbf{e}_s$.
Here we impose the canopy node condition, $\hat{\mathbf{c}}={\mathbf{c}}$ or $A\mathbf{1}=\widehat{A}\mathbf{1}$,
such that the IERK method \eqref{Scheme: general IERK} is consistent at all stages. Table \ref{table: order condition} lists the order conditions for the coefficients matrices and the weight vectors to make the IERK method \eqref{Scheme: general IERK}  to be accurate up to third-order in time. A detailed description of these order conditions can also be found in
\cite{AscherRuuthSpiteri:1997,IzzoJackiewicz:2017,ShinLeeLee:2017}.

\begin{table}[htb!]
	\centering
	\begin{threeparttable}
		\centering 
		\renewcommand\arraystretch{1.3}
		\belowrulesep=0pt\aboverulesep=0pt
		\caption{Order conditions for IERK methods up to third-order.}
		\label{table: order condition}
		\vspace{2mm}
		\begin{tabular}{c|cc|c}
			\toprule 
			\multirow{2}*{Order} & \multicolumn{2}{c|}{Stand-alone conditions} & \multirow{2}*{Coupling condition} \\
			\cmidrule{2-3}
			& Implicit part & Explicit part & \\
			\midrule 
			1 & $\mathbf{b}^{T} \mathbf{1}=1$ & $\hat{\mathbf{b}}^{T} \mathbf{1} = 1$ & - \\[2pt]
			\hline 2 & $\mathbf{b}^{T} \mathbf{c} = \tfrac1{2}$ & $\hat{\mathbf{b}}^{T} \mathbf{c} = \tfrac1{2}$ & - \\[2pt]
			\hline 3 & $\mathbf{b}^{T} \mathbf{c}^{.2} = \tfrac1{3}$ & $\hat{\mathbf{b}}^{T} \mathbf{c}^{.2} = \tfrac1{3}$ & \\
			& $\mathbf{b}^{T} A  \mathbf{c} = \tfrac1{6}$ & $\hat{\mathbf{b}}^{T} \widehat{A}  \mathbf{c} =\tfrac1{6}$ & {$\mathbf{b}^{T} \widehat{A}  \mathbf{c} = \tfrac1{6}, \; \hat{\mathbf{b}}^{T} A  \mathbf{c} = \tfrac1{6}$} \\[2pt]
			\bottomrule
		\end{tabular}
		\footnotesize
		\tnote{* For the vectors $\mathbf{x}$ and $\mathbf{y}$, $\mathbf{x}\odot\mathbf{y}:=(x_1y_1,x_2y_2,\cdots,x_sy_s)^T$ and $\mathbf{x}^{.m}:=\mathbf{x}\odot\mathbf{x}^{.(m-1)}$ for $m>1$.}
	\end{threeparttable}
\end{table}

Moreover, requiring $u_h^{n,i}=\phi_h^*$ for all
$i$ and $n \geq 1$ immediately shows that the canopy node condition $\hat{\mathbf{c}}={\mathbf{c}}$ makes the IERK method \eqref{Scheme: general IERK} preserve naturally the equilibria $\phi^*$ of the CH model \eqref{cont: Problem-CH stabilized}, that is, $\epsilon^2\Delta_h\phi_h^*=F'(\phi_h^*)$ or $L_{\kappa,h}\phi_h^*=f_{\kappa}(\phi_h^*)$. So one can reformulate the standard form \eqref{Scheme: general IERK} into the following steady-state preserving form
\begin{align}\label{Scheme: IERK steady-state preserving}				&\,u_h^{n,i+1}=u_h^{n,1}+\tau_n\sum_{j=1}^{i+1}a_{i+1,j}\Delta_h \brab{L_{\kappa,h}u_h^{n,j}-L_{\kappa,h}u_h^{n,1}}-\tau_n\sum_{j=1}^{i}\hat{a}_{i+1,j}\Delta_h \kbrab{f_{\kappa}(u_h^{n,j})-L_{\kappa,h}u_h^{n,1}}\nonumber\\
	&\,=u_h^{n,1}+\tau_n\sum_{j=1}^{i}a_{i+1,j+1}\Delta_h \brab{L_{\kappa,h}u_h^{n,j+1}-L_{\kappa,h}u_h^{n,1}}-\tau_n\sum_{j=1}^{i}\hat{a}_{i+1,j}\Delta_h \kbrab{f_{\kappa}(u_h^{n,j})-L_{\kappa,h}u_h^{n,1}}
\end{align}
for $1\le i\le s_{\mathrm{I}}:=s-1$ ($s_{\mathrm{I}}$ represents the number of implicit stages), in which we drop the terms with the coefficients $a_{i+1,1}$ for $1\le i\le s_{\mathrm{I}}$.
In this sense, we define the lower triangular coefficient matrices for the implicit and explicit parts, respectively,
$$A_{\mathrm{I}}:=\brab{a_{i+1,j+1}}_{i,j=1}^{s_{\mathrm{I}}}\quad\text{and}\quad A_{\mathrm{E}}:=\brab{\hat{a}_{i+1,j}}_{i,j=1}^{s_{\mathrm{I}}}\,.$$
Note that, the two matrices $A_{\mathrm{I}}$ and $A_{\mathrm{E}}$ are always required in our theory, while the coefficient vector $\mathbf{a}_{1}:=(a_{21},a_{31},\cdots,a_{s1})^T$ would be not involved directly although it would be useful in designing some computationally effective IERK methods. If $\mathbf{a}_{1}\neq\mathbf{0}$, the Lobatto-type DIRK methods are used in the implicit part and the associated method \eqref{Scheme: general IERK} is called Lobatto-type IERK methods; while we call \eqref{Scheme: general IERK} as Radau-type  (known as ARS-type \cite{AscherRuuthSpiteri:1997}) IERK  methods if $\mathbf{a}_{1}=\mathbf{0}$.

Actually, to avoid nonlinear iteration at each stage of implicit Runge-Kutta methods, the IERK methods have attracted much attention \cite{AscherRuuthSpiteri:1997,CardoneJackiewiczSanduZhang:2014MMA,FuTangYang:2024,IzzoJackiewicz:2017,KennedyCarpenter:2003,LiaoWangWen:2024arxiv,ShinLeeLee:2017,ShinLeeLee:2017CMA}.
Kennedy and Carpenter \cite{KennedyCarpenter:2003} constructed high-order IERK methods from third- to fifth-order to simulate convection-diffusion-reaction equations. The widespread ARS-type IERK methods were developed in \cite{AscherRuuthSpiteri:1997} to solve the convection-diffusion problems.  
Cardone et al. \cite{CardoneJackiewiczSanduZhang:2014MMA} proposed a class of ARS-type IERK methods up to fourth-order based on the extrapolation of stage solutions at the current and previous steps. Izzo and Jackiewicz \cite{IzzoJackiewicz:2017} constructed some parameterized IERK methods up to fourth-order with A-stable implicit part by choosing the method parameters to maximize the regions of absolute stability for the explicit part. In simulating the semilinear parabolic problems, IERK methods turned out to be very competitive. Shin et al. \cite{ShinLeeLee:2017CMA} observed that the Radau-type IERK methods combined with convex splitting technique exhibit the original energy stability in numerical experiments, and established the energy stability in \cite{ShinLeeLee:2017} for a special case. 
Recently, Fu et al. \cite{FuTangYang:2024} derived some sufficient conditions of Radau-type IERK methods to maintain the decay of original energy for the gradient flow problems with the global Lipschitz continuity assumption of nonlinear bulk, and presented some concrete schemes up to third-order accuracy. However, the rigorous stability and convergence of these (not algebraically stable) IERK methods for the CH model \eqref{cont: Problem-CH} are rather difficult because it is technically challenging to verify the uniform boundedness of stage solutions for multi-stage methods \cite{FuTangYang:2024,ShinLeeLee:2017,ShinLeeLee:2017CMA}.

In the literature, there are many works on linearized time-stepping methods and the numerical analysis for the CH model \eqref{cont: Problem-CH}. Sun \cite{Sun:1995} developed three-level Crank-Nicolson type method with finite difference approximation in space, and Li et al. \cite{LiSunZhao:2012} established the unconditional convergence of the three-level linearized time-stepping method in the maximum norm. Second-order energy-stable semi-implicit approaches were constructed and analyzed in \cite{HeLiuTang:2007StabilizedCH,WangYu:2018OnStabilizedCH}.  In resolving the long-time multi-scale dynamics, variable-step time-stepping approaches were developed and analyzed in \cite{ChenWangYanZhang:2019,LiaoJiWangZhang:2022,ZhangQiao:2012}, in which the robust stability and convergence with respect to the change of time-step sizes are established. Combining the scalar auxiliary variable method \cite{ShenXuYang:2018} and the Gaussian Runge-Kutta methods, arbitrarily high-order energy-stable schemes were constructed in \cite{AkrivisLiLi:2019}, while the energy stability was established with respect to a modified energy involving the auxiliary variable. Related works on energy stable time-stepping methods can be also found in \cite{ChengWangWiseYue:2016Weakly,FuTangYang:2024,HuangYangWei:2020,LiQiaoWang:2021,ShenXuYang:2018} and the reference therein, with certain modified discrete energy functional by adding some nonnegative small terms to the original energy \eqref{cont:free energy}, especially for high-order multi-step approaches.

In the next section, we present the theoretical and numerical motivations to a refined class of Lobatto-type IERK methods, called refined IERK (R-IERK), in which the associated differentiation matrices and the average dissipation rates are independent of the time-space discretization meshes. Some parameterized R-IERK methods up to third-order accuracy are constructed in Section \ref{sec: Construction R-IERK}. The unconditional $L^2$ norm error estimate and the original energy dissipation law of R-IERK methods are addressed in Section \ref{sec: convergence} with an updated technique of time-space error splitting. 
Extensive experiments are presented in Section \ref{sec: experiments} to support our theory and show the effectiveness of R-IERK methods. 

For simplicity, we will use the simplified notations, IERK($p$, $s$; $\mu$) and R-IERK($p$, $s$; $\mu$), to represent the  $p$-th order $s$-stage $\mu$-parameterized IERK and R-IERK methods, respectively. Always, we use $D_{\mathrm{R}}^{(p,s)}$ and $\mathcal{R}_{\mathrm{R}}^{(p,s)}$ to represent the associated differentiation matrix and the average dissipation rate, respectively, of a  $p$-th order $s$-stage R-IERK method. 


\section{Theoretical and numerical motivations}\label{sec: motivations}
\setcounter{equation}{0}


Set the domain $\Omega=(0,L)^2$
and consider the uniform length $h_x=h_y=h:=L/M$ in each direction
for an even positive integer $M$.
Let $\Omega_{h}:=\big\{\myvec{x}_{h}=(ih,jh)\,|\,1\le i,j \le M\big\}$
and put
$\bar{\Omega}_{h}:=\Omega_{h}\cup\partial{\Omega}$.
Denote the space of $L$-periodic grid functions
$\mathbb{V}_{h}:=\{v\,|\,v=\bra{v_h}\; \text{is $L$-periodic for}\; \myvec{x}_h\in\bar{\Omega}_h\}.$
For a periodic function $v(\myvec{x})$ on $\bar{\Omega}$,
let $P_M:L^2(\Omega)\rightarrow \mathscr{F}_M$
be the standard $L^2$ projection operator onto the space $\mathscr{F}_M$,
consisting of all trigonometric polynomials of degree up to $M/2$,
and $I_M:L^2(\Omega)\rightarrow \mathscr{F}_M$
be the trigonometric interpolation operator \cite{ShenTangWang:2011Spectral},
i.e.,
\[
\bra{P_Mv}(\myvec{x})=\sum_{\ell,m =- M/2}^{M/2-1}
\widehat{v}_{\ell,m}e_{\ell,m}(\myvec{x}),\quad
\bra{I_Mv}(\myvec{x})=\sum_{\ell,m=- M/2}^{M/2-1}
\widetilde{v}_{\ell,m}e_{\ell,m}(\myvec{x}),
\]
where the complex exponential basis function
$e_{\ell,m}(\myvec{x}):=e^{\mathrm{i}\nu\bra{\ell x+my}}$ with $\nu=2\pi/L$.
The coefficients $\widehat{v}_{\ell,m}$
refer to the standard Fourier coefficients of function $v(\myvec{x})$,
and the
pseudo-spectral coefficients $\widetilde{v}_{\ell,m}$ are determined such that $\bra{I_Mv}(\myvec{x}_h)=v_h$. 
The Fourier pseudo-spectral first and second order derivatives of $v_h$ are given by
\[
\mathcal{D}_xv_h:=\sum_{\ell,m= -M/2}^{M/2-1}
\bra{\mathrm{i}\nu\ell}\widetilde{v}_{\ell,m}
e_{\ell,m}(\myvec{x}_h),\quad
\mathcal{D}_x^2v_h:=\sum_{\ell,m = -M/2}^{M/2-1}
\bra{\mathrm{i}\nu\ell}^2\widetilde{v}_{\ell,m}
e_{\ell,m}(\myvec{x}_h).
\]
The operators $\mathcal{D}_y$
and $\mathcal{D}_y^2$ can be defined in the similar fashion.
We define the discrete gradient and Laplacian
in the point-wise sense by $\nabla_hv_h := \left(\mathcal{D}_xv_h,\mathcal{D}_yv_h\right)^T$ and $\Delta_hv_h :=\bra{\mathcal{D}_x^2+\mathcal{D}_y^2}v_h,$ respectively.

For functions $v,w\in\mathbb{V}_{h}$,
define the discrete inner product
$\myinner{v,w}:=h^2\sum_{\myvec{x}_h\in\Omega_{h}}v_h w_h$,
and the associated discrete $L^{2}$ norm $\mynorm{v}:=\mynorm{v}_{l^2}=\sqrt{\myinner{v,v}}$.
Also, we will use  $\mynorm{v}_{\infty}=\max_{\myvec{x}_h\in\Omega_{h}}|v_h|$,
$\mynormb{\nabla_hv}:=\sqrt{h^2\sum_{\myvec{x}_h\in\Omega_{h}}|\nabla_hv_h|^2}$
 and  $\mynormb{\Delta_hv}:=\sqrt{h^2\sum_{\myvec{x}_h\in\Omega_{h}}|\Delta_hv_h|^2}$.
It is easy to check the discrete Green's formulas,
$\myinner{-\Delta_hv,w}=\myinner{\nabla_hv,\nabla_hw}$ and
$\myinnert{\Delta_h^2v,w}=\myinner{\Delta_hv,\Delta_hw}$,
see \cite{ShenTangWang:2011Spectral,ChengWangWiseYue:2016Weakly} for more details.
Also one has the embedding inequality simulating the Sobolev embedding $H^2(\Omega)\hookrightarrow L^\infty(\Omega)$,
\begin{align}\label{ieq: H2 embedding}
	\mynormb{v}_{\infty}\le \ck_\Omega\brab{\mynormb{v}+\mynormb{\Delta_hv}}\quad \text{for any $v\in\mathbb{V}_{h}$,}
\end{align}
where $\ck_{\Omega}>0$ is only dependent on the domain $\Omega_h$. Here and hereafter, any subscripted $\ck$, such as $\ck_\Omega, \ck_u, \ck_\phi$, $\ck_0$, $\ck_1$, $\hck_1$, $\hck_1^*$ and so on,
denotes a fixed constant; while any subscripted $\Ck$, such as $\Ck_w$ and $\Ck_\phi$, denotes a generic positive constant, not necessarily
the same at different occurrences. The appeared constants may be dependent on the given data
(typically, the interface width parameter $\epsilon$)
and the solution $\Phi$ but are always independent of the spatial length $h$ and the time-step size $\tau_n$.

\subsection{Theoretical motivation}

Let $E_{s_{\mathrm{I}}}:=(1_{i\ge j})_{s_{\mathrm{I}}\times s_{\mathrm{I}}}$ be the lower triangular matrix full of element 1, and $I_{s_{\mathrm{I}}}$ be the  identity matrix of the same size as $A_{\mathrm{I}}$. Following the derivations in \cite[Section 2]{LiaoWangWen:2024arxiv}, one can reformulate the IERK method \eqref{Scheme: IERK steady-state preserving} into the following differential form,
also see \cite{LiaoWang:2024MCOM,LiaoWangWen:2024JCP},
\begin{align}\label{Scheme: DOC IERK stabilized}			
	\sum_{\ell=1}^{i}d_{i\ell}\brab{\tau_n \Delta_h  L_{\kappa,h}}\delta_{\tau}u_h^{n,\ell+1}    
	=\tau_n \Delta_h \kbraB{L_{\kappa,h}u_h^{n,i+\frac12}-f_{\kappa}(u_h^{n,i})}\quad\text{for $n\ge1$ and $1\le i\le s_{\mathrm{I}}$,}
\end{align}
where $\delta_{\tau}u_h^{n,\ell+1}:=u_h^{n,\ell+1}-u_h^{n,\ell}$ for $\ell\ge 1$  is the time difference, $u_h^{n,\ell+\frac12}:=(u_h^{n,\ell+1}+u_h^{n,\ell})/2$ for $\ell\ge 1$  is the averaged operator  and $d_{i\ell}$ is the element of the differentiation matrix  $D=(d_{i\ell})_{s_{\mathrm{I}}\times s_{\mathrm{I}}}$ defined by
\begin{align}\label{Def: Differential Matrix D}			
	D(z):=D_{\mathrm{E}}-zD_{\mathrm{EI}}\quad\text{with}\quad D_{\mathrm{E}}:=A_{\mathrm{E}}^{-1}E_{s_{\mathrm{I}}},\quad D_{\mathrm{EI}}:=A_{\mathrm{E}}^{-1}A_{\mathrm{I}}E_{s_{\mathrm{I}}}
	-E_{s_{\mathrm{I}}}+\tfrac{1}2I_{s_{\mathrm{I}}}.
\end{align}

It is easy to check that the stage solutions $u_h^{n,i+1}$ of the IERK method \eqref{Scheme: IERK steady-state preserving} or \eqref{Scheme: DOC IERK stabilized} preserve the initial volume, that is, 
\begin{align}\label{eq: stage volume preserving}	
\myinnerb{u^{n,i+1},1}=\myinnerb{u^{n,1},1}=\myinnerb{u^{1,1},1}=\myinnerb{\phi^{0},1}
\quad\text{for $n\ge1$ and $1\le i\le s_{\mathrm{I}}$}.
\end{align}
To establish the original energy dissipation law of the IERK method \eqref{Scheme: DOC IERK stabilized}, we assume that the nonlinear term $F'(\Phi)$ is continuously differentiable and  recall the following inequality 
\begin{align*}	
	F'(u_h^{n,i})(u_h^{n,i+1}-u_h^{n,i})
	\le &\,F(u_h^{n,i})-F(u_h^{n,i+1})
	+\tfrac1{2}(u_h^{n,i+1}-u_h^{n,i})^2\max_{\xi_h\in \mathcal{B}_{n,i}}\abs{F''(\xi_h)}
\end{align*}
where the function space $\mathcal{B}_{n,i}:=\big\{\xi_h\,\big|\big.\,\mynormb{\xi}_{\infty}\le \max\{\mynormt{u^{n,i+1}}_{\infty},\mynormt{u^{n,i}}_{\infty}\}\big\}$. By using the definitions in \eqref{def: stabilized parameter}, it is easy to find that
\begin{align*}
	\myinnerb{L_{\kappa,h}u^{n,i+\frac12}-f_{\kappa}(u^{n,i}),\delta_{\tau}u^{n,i+1}}
	\le E[u^{n,i}]-E[u^{n,i+1}]
	-\tfrac1{2}\mynormb{\delta_{\tau}u^{n,i+1}}^2
	\brab{\kappa-\max_{\xi_h\in \mathcal{B}_{n,i}}\mynormb{F''(\xi)}_{\infty}}.
\end{align*}
Then by imposing the uniformly bounded assumption of stage solutions and taking the stabilized parameter $\kappa\ge\max_{\xi_h\in \mathcal{B}_{n,i}}\mynormb{F''(\xi)}_{\infty}$ (the minimum size of $\kappa$ would be always dependent on the regularity of initial data and the magnitude of the small interface parameter $\epsilon$, see \cite{LiQiao:CMS2017,LiQiao:JSC2017,LiQiaoTang:SINUM2018,LiQuanTang:MC2022}), one can follow the proofs of \cite[Theorem 2.1, Corollary 2.1 and Lemma 2.1]{LiaoWangWen:2024arxiv} to obtain the following result. Here and hereafter, we say that a lower triangular matrix $D$ is positive (semi-)definite
if its symmetric part $\mathcal{S}(D)=(D+D^T)/2$ is  positive (semi-)definite.

\begin{lemma}\label{lemma: IERK energy stability} 
	Assume that the two matrices $D_{\mathrm{E}}$ and $D_{\mathrm{EI}}$ in \eqref{Def: Differential Matrix D}	 are positive (semi-)definite. If the stage solutions $u_h^{n,\ell}$ $(n\ge1,1\le\ell\le s)$ are bounded by $\ck_0$ in the maximum norm, and the stabilization parameter $\kappa$ in \eqref{def: stabilized parameter} is chosen properly large such that 
	$\kappa\ge\max_{\mynormt{\xi}_{\infty}\le \ck_0}\mynormb{F''(\xi)}_{\infty},$
	then the IERK method \eqref{Scheme: DOC IERK stabilized} preserves 
	the original energy dissipation law \eqref{cont:energy dissipation}
	at all stages,
	\begin{align}\label{thmResult: stage energy laws}			
		E[u^{n,i+1}]-E[u^{n,1}]\le&\,
		\frac1{\tau_n}\sum_{k=1}^{i}\myinnerB{\Delta_h^{-1}\delta_{\tau}u^{n,k+1},
			\sum_{\ell=1}^{k}d_{k\ell}({\tau_n}\Delta_h L_{\kappa,h})\delta_{\tau}u^{n,\ell+1}}
	\end{align}
	for $n\ge1$ and $1\le i\le s_{\mathrm{I}}$.
	The associated average dissipation rate is nonnegative, that is,
	\begin{align}\label{def: numerical rate Formula}		
		\mathcal{R}=\frac{1}{s_{\mathrm{I}}}\mathrm{tr}(D_{\mathrm{E}})+\frac{1}{s_{\mathrm{I}}}\mathrm{tr}(D_{\mathrm{EI}})\tau_n\overline{\lambda}_{\mathrm{ML}}=\frac1{s_{\mathrm{I}}}\sum_{k=1}^{s_{\mathrm{I}}}\frac1{\hat{a}_{k+1,k}}
		+\frac{1}{s_{\mathrm{I}}}\sum_{k=1}^{s_{\mathrm{I}}}
		\braB{\frac{a_{k+1,k+1}}{\hat{a}_{k+1,k}}-\frac12}\tau_n\overline{\lambda}_{\mathrm{ML}}\ge0,  
	\end{align}   
	where $\overline{\lambda}_{\mathrm{ML}}>0$ is the average eigenvalue of the symmetric, positive definite matrix $-\Delta_hL_{\kappa,h}$.    	
\end{lemma}

By comparing the discrete energy law \eqref{thmResult: stage energy laws} with the continuous counterpart \eqref{cont:energy dissipation}, as remarked in \cite{LiaoWang:2024MCOM,LiaoWangWen:2024arxiv},  an IERK method would be a ``good" candidate if the average dissipation rate $\mathcal{R}$ is as close to 1 as possible within properly large range of $\tau_n\overline{\lambda}_{\mathrm{ML}}$.
Always, potential users will determine the value of $\tau_n\overline{\lambda}_{\mathrm{ML}}$ by choosing the time-step size $\tau_n$, the stabilized parameter $\kappa$, the method of spatial approximation and the spacing size $h$ as well. \textit{This raises an interesting question: are there any IERK methods whose average dissipation rates are independent of $\tau_n\overline{\lambda}_{\mathrm{ML}}$?}

It has a positive answer which gives rise to a refined class of IERK method \eqref{Scheme: DOC IERK stabilized}, called refined IERK (R-IERK), in which the associated  average dissipation rate $\mathcal{R}$ in \eqref{def: numerical rate Formula} is always independent of $\tau_n\overline{\lambda}_{\mathrm{ML}}$.  Actually, from the formula  \eqref{def: numerical rate Formula}, we see that $\mathcal{R}$ is independent of $\tau_n\overline{\lambda}_{\mathrm{ML}}$ if and only if the positive semi-definite matrix $D_{\mathrm{EI}}$ is a zero matrix, that is,
\begin{align}\label{setting: zero Dei}			
	D_{\mathrm{EI}}=A_{\mathrm{E}}^{-1}A_{\mathrm{I}}E_{s_{\mathrm{I}}}
	-E_{s_{\mathrm{I}}}+\tfrac{1}2I_{s_{\mathrm{I}}}=\mathbf{0}.
\end{align}
It can be achieved by taking $A_{\mathrm{I}}:=A_{\mathrm{E}}P_{s_{\mathrm{I}}}$ with
$P_{s_{\mathrm{I}}}:=I_{s_{\mathrm{I}}}-\tfrac{1}2E_{s_{\mathrm{I}}}^{-1}.$ 
In such case,  the differentiation  matrix $D(z)$ in \eqref{Def: Differential Matrix D} is independent of $z$,
that is, 
\begin{align}\label{def: DOC dij}		
	D_{\mathrm{R}}:=\brab{d_{ij}^{(R)}}_{s_{\mathrm{I}}\times s_{\mathrm{I}}}=D_{\mathrm{E}}=A_{\mathrm{E}}^{-1}E_{s_{\mathrm{I}}},
\end{align} 
and the associated average dissipation rate $\mathcal{R}$ in \eqref{def: numerical rate Formula} is reduced into 
\begin{align}\label{def: numerical rate R-IERK}		
	\mathcal{R}_{\mathrm{R}}:=\frac1{s_{\mathrm{I}}}
	\sum_{k=1}^{s_{\mathrm{I}}}\frac1{\hat{a}_{k+1,k}}.  
\end{align} 
The IERK scheme \eqref{Scheme: IERK steady-state preserving} reduces into 
the following R-IERK method
\begin{align}\label{Scheme: R-IERK stabilized}				
	u_h^{n,i+1}=&\,u_h^{n,1}+\tau_n\sum_{j=1}^{i}\hat{a}_{i+1,j}\Delta_h \kbraB{L_{\kappa,h}u_h^{n,j+\frac1{2}}-f_{\kappa}(u_h^{n,j})}
	\quad\text{for $n\ge1$ and $1\le i\le s_{\mathrm{I}}$,}
\end{align}
or, with the matrix $(\underline{\hat{a}}_{i+1,j})_{s_{\mathrm{I}}\times s_{\mathrm{I}}}:=E_{s_{\mathrm{I}}}^{-1}A_{\mathrm{E}}$,
\begin{align}\label{Scheme: R-IERK difference}				
	\delta_{\tau}u_h^{n,i+1} 
	=&\,\tau_n\sum_{j=1}^{i}\underline{\hat{a}}_{i+1,j}\Delta_h \kbraB{L_{\kappa,h}u_h^{n,j+\frac1{2}}-f_{\kappa}(u_h^{n,j})}
	\quad\text{for $n\ge1$ and $1\le i\le s_{\mathrm{I}}$.}
\end{align}

	
The implicit part of the IERK method \eqref{Scheme: DOC IERK stabilized} will be completely determined by its explicit part. It means that more stages and computational cost would be required to achieve the interesting property \eqref{setting: zero Dei} of R-IERK methods due to the vanish of the degree of freedom in the implicit part. According to the canopy node condition, $\hat{\mathbf{c}}={\mathbf{c}}$ or $A\mathbf{1}=\widehat{A}\mathbf{1}$, it is not difficult to compute the coefficient vector of the first column, $\mathbf{a}_{1}=\brab{\tfrac12\hat{a}_{21},\tfrac12\hat{a}_{31},\cdots,\tfrac12\hat{a}_{s1}}^T.$
As seen, the implicit part of the R-IERK method \eqref{Scheme: R-IERK stabilized} allows a stage-order of two so that it would be	 not necessarily algebraically stable, cf.  \cite{KennedyCarpenter:2003,Spijker:1985}. Moreover, the non-zero vector $\mathbf{a}_{1}$ implies the following result because the Radau-type or ARS-type IERK methods always require $\mathbf{a}_{1}=\mathbf{0}$, see \cite{AscherRuuthSpiteri:1997,CardoneJackiewiczSanduZhang:2014MMA,FuTangYang:2024}.
\begin{proposition}
	There are no Radau-type or ARS-type IERK methods that the associated average dissipation rates are independent of the time-space discretization parameter $\tau_n\overline{\lambda}_{\mathrm{ML}}$.
\end{proposition}   	

\subsection{Numerical motivation}\label{subsec: motivation}

Before we turn to high-order R-IERK methods, 
some preliminary tests are presented to show the advantage of R-IERK methods in self-adaptive computations. Consider the IERK(1,2;$\theta$) method
\begin{equation}\label{Scheme: IERK1}
	\begin{array}{c|c}
		\mathbf{c} & A \\
		\hline\\[-10pt]   & \mathbf{b}^T
	\end{array}
	=\begin{array}{c|cc}
		0 & 0 &  0  \\
		1 & 1-\theta & \theta \\
		\hline  & 1-\theta & \theta
	\end{array}\;,\quad \begin{array}{c|c}
		\hat{\mathbf{c}} & \widehat{A} \\
		\hline\\[-10pt]   & \hat{\mathbf{b}}^T
	\end{array}
	=\begin{array}{c|cc}
		0 & 0 & 0   \\
		1 & 1 & 0 \\
		\hline  & 1 & 0 
	\end{array}\;.
\end{equation}
Obviously, $A_{\mathrm{I}}=(\theta)$ and $A_{\mathrm{E}}=(1)$ such that the differentiation matrix $D^{(1,2)}(z)=1-z\bra{\theta-\tfrac{1}2}$ is positive definite for $z\le0$ provided  $\theta\ge\frac12$ according to Lemma \ref{lemma: IERK energy stability}. The average dissipation rate
\begin{align}\label{AverRate: IERK1}
	\mathcal{R}^{(1,2)}(\theta)=1+\bra{\theta-\tfrac{1}2}\tau_n\overline{\lambda}_{\mathrm{ML}}\quad\text{for $\theta\ge\tfrac12$}.
\end{align}
Choosing the weighted parameter $\theta=1/2$, we have $\mathcal{R}^{(1,2)}_{\mathrm{R}}=\mathcal{R}^{(1,2)}(\tfrac12)=1$, which is independent of  $\tau_n\overline{\lambda}_{\mathrm{ML}}$, and the associated Crank-Nicolson-type R-IERK(1,2) scheme,
\begin{align*}
	\delta_{\tau}\phi_h^{n}    
	=\tau_n \Delta_h\kbra{\tfrac{1}2L_{\kappa,h}(\phi_h^{n}+\phi_h^{n-1})-f_{\kappa}(\phi_h^{n-1})}\quad \text{for $n\ge1$.}  
\end{align*}

\begin{example}\label{1example: energy_2}
	Consider the Cahn-Hilliard model \eqref{cont: Problem-CH}
	on the spatial domain $\Omega=(-\pi,\pi)^2$ with the interface parameter $\epsilon=0.1$, subject to the initial data 
	\begin{align*}
		\Phi^{0} =&\, \tfrac{1}{2} \tanh(|x| + |y| + 1) - e^{-5(|x| + |y| - 2)^2} + \tfrac{1}{2} e^{-2(|x| - 1)^2} + \tfrac{1}{10}\sin(e^{|y| - 1}).
	\end{align*}
	The reference solution is generated with $\tau=10^{-4}$ by the Lobatto-type IERK2-2 method in \cite{LiaoWangWen:2024arxiv} for the parameter $a_{33}=\tfrac{1+\sqrt{2}}{4}$, called IERK(2,3) method hereafter, which is regarded as the best one among second-order IERK methods in \cite{LiaoWangWen:2024arxiv} with the associated average dissipation rate $\mathcal{R}^{(2,3)}_{\mathrm{L}}
	=\sqrt{2}+\tfrac{\sqrt{2}}{4}\tau_n\overline{\lambda}_{\mathrm{ML}}$.
\end{example}

To investigate the energy behavior in adaptive time-stepping calculation, we always choose the adaptive step-size \cite{HuangYangWei:2020,LiaoJiWangZhang:2022}, $\tau_{ada}=\max\{\tau_{\min},\tau_{\max}/\Pi_\eta(\phi)\}$ with $\Pi_\eta(\phi):=\sqrt{1+\eta\|\partial_\tau\phi^n\|^2}$, where	\( \tau_{\max} \) and \( \tau_{\min} \) refer to the predetermined maximum and minimum time-step sizes, and $\eta$ is chosen by the user to adjust the level of adaptivity. Hereafter, if not explicitly specified, we set the adaptivity parameter $\eta=1000$, the minimum step $\tau_{\min}=10^{-4}$ and the initial step size $\tau_1=\tau_{\min}$.


\begin{figure}[htb!]
	\centering
	\subfigure[$\tau_{\max}=0.01$]{
		\includegraphics[width=2.05in]{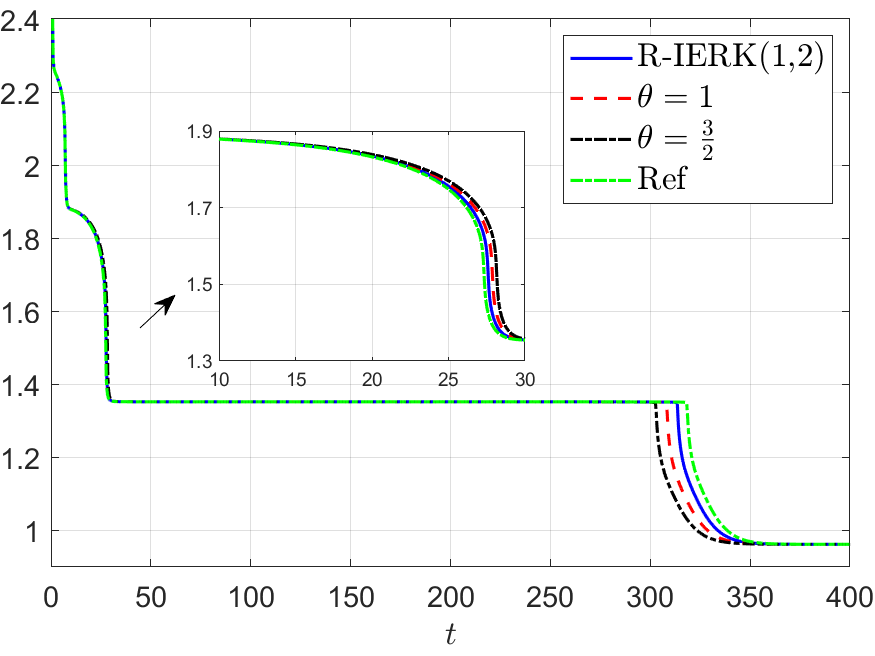}}
	\subfigure[$\tau_{\max}=0.05$]{
		\includegraphics[width=2.05in]{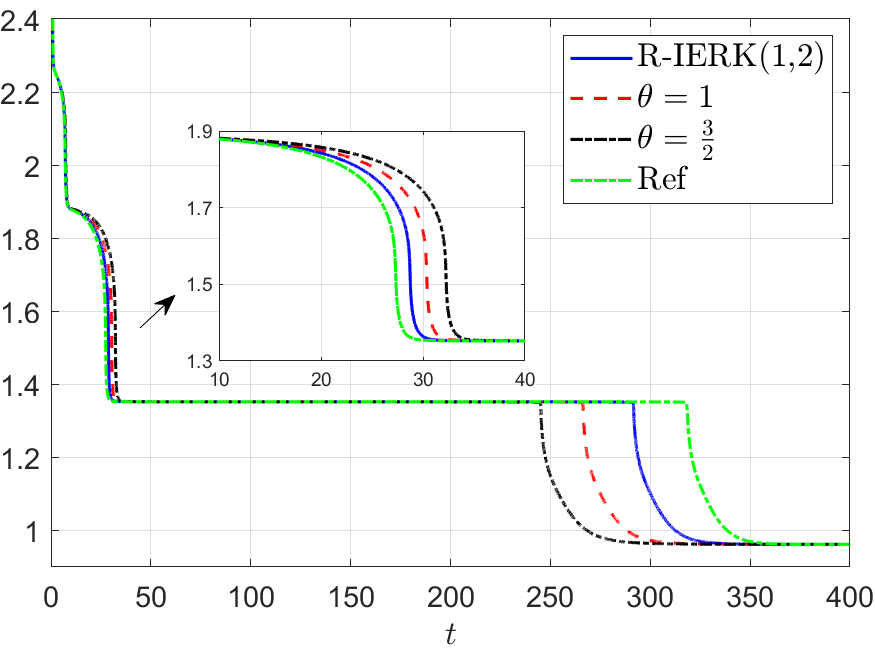}}
	\subfigure[$\tau_{\max}=0.1$]{
		\includegraphics[width=2.05in]{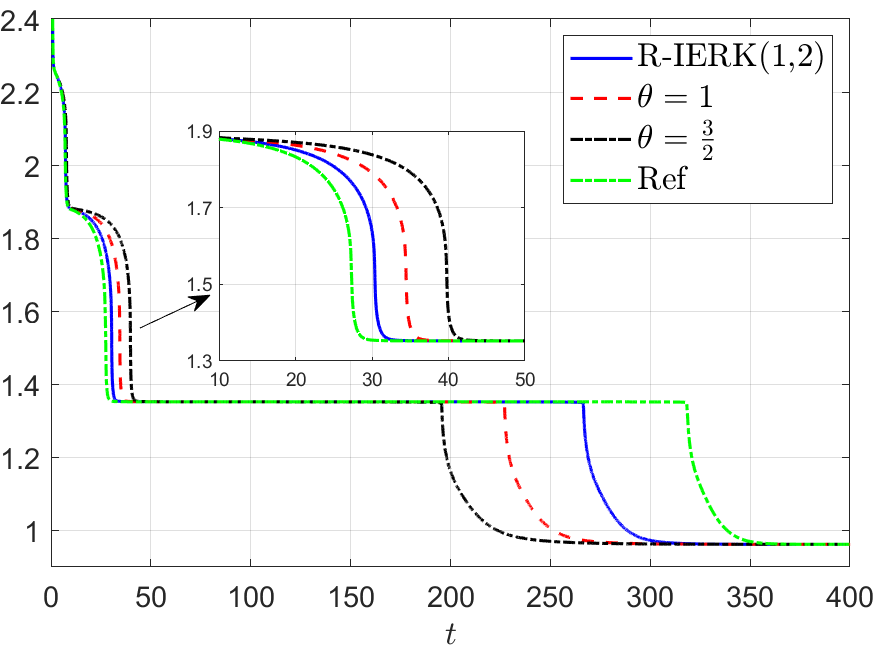}}
	\caption{Energy behaviors of R-IERK(1,2) method and IERK(1,2;$\theta$) method \eqref{Scheme: IERK1}.}
	\label{fig: IERK1 decay theta}
\end{figure}

We run the R-IERK(1,2) method and the IERK(1,2;$\theta$) method \eqref{Scheme: IERK1} with the stabilized parameter $\kappa=2.5$ to the final time $T=1000$. Figure \ref{fig: IERK1 decay theta} depicts the discrete energy curves (the discrete energies over $t\in[400,1000]$ are unchanged and omitted here to display the rapidly changing parts more finely) for three different scenes: (a) $\tau_{\max}=0.01$, (b) $\tau_{\max}=0.05$ and (c) $\tau_{\max}=0.1$.
Under the same time-adaptive strategy with different maximum time step $\tau_{\max}$, 
the energy curves of the R-IERK(1,2) method are always closer to the reference energy than those of the IERK(1,2;$\theta$) methods with $\theta=1$ and $\frac32$. The IERK(1,2;$\theta$) methods with $\theta=1$ and $\frac32$ seem more susceptible with respect to the change of time-step sizes; while the R-IERK(1,2) method allows some larger adaptive step-sizes and would be more preferable in adaptive numerical simulations. This phenomenon would be attributed to the constant rate of average dissipation rate, $\mathcal{R}^{(1,2)}_{\mathrm{R}}=1$, which is independent of the time-space discretization parameter $\tau_n\overline{\lambda}_{\mathrm{ML}}$.
  
  \section{Construction of parameterized R-IERK methods}
  \label{sec: Construction R-IERK}
  \setcounter{equation}{0}

 \subsection{Second-order R-IERK methods}

 Second-order R-IERK methods require two implicit stages, $s_I=2$, at least; however, it is easy to check that no 3-stage R-IERK method exists for second-order accuracy.
 
 We turn to consider 4-stage R-IERK methods that satisfy the canopy node condition and the first-order conditions in Table \ref{table: order condition},
 \begin{equation*}
 	\begin{array}{c|c}
 		\mathbf{c} & A \\
 		\hline\\[-10pt] 	& \mathbf{b}^T
 	\end{array}
 	=	\begin{array}{c|cccc}
 		0 & 0 &  &    &  \\[1pt]
 		c_{2} & \tfrac{c_{2}}{2} & \tfrac{c_{2}}{2} &   \\[1pt]
 		c_3 & \tfrac{c_3-c_2}{2} & \tfrac{c_3}{2} & \tfrac{c_2}{2}    \\[1pt]
 		1 & \frac{1-\hat{a}_{42}-c_2}{2} & \tfrac{1-c_2 }{2} & \tfrac{\hat{a}_{42}+c_2 }{2} &   \tfrac{c_2}{2} \\[1pt]
 		\hline\\[-10pt]   & \frac{1-\hat{a}_{42}-c_2}{2} & \tfrac{1-c_2 }{2} & \tfrac{\hat{a}_{42}+c_2 }{2} &   \tfrac{c_2}{2}
 	\end{array},\quad
 	\begin{array}{c|c}
 		\hat{\mathbf{c}} & \widehat{A} \\
 		\hline\\[-10pt] 	& \hat{\mathbf{b}}^T
 	\end{array}
 	=\begin{array}{c|cccc}
 		0 & 0 &  &   &  \\[1pt]
 		c_{2} & c_2 & 0 & &  \\[1pt]
 		c_3 & c_3-c_2 &c_2 & 0&   \\[1pt]
 		1 & 1-\hat{a}_{42}-c_2 &\hat{a}_{42} & c_2 & 0 \\[1pt]
 		\hline\\[-10pt]   & 1-\hat{a}_{42}-c_2 &\hat{a}_{42} & c_2&0
 	\end{array},
 \end{equation*}
 where we set $\hat{a}_{32}:=c_2$ and $\hat{a}_{43}:=c_2$ to reduce the degree of freedom
 and there are three independent coefficients $c_2$, $c_3$ and $\hat{a}_{42}$. These simplifying settings are also useful for practical applications since they provide a constant iteration matrix for the system of linear equations at all stages. According to the stand-alone conditions for second-order accuracy, $\hat{\mathbf{b}}^T\mathbf{c}=\tfrac12$ and $\mathbf{b}^T\mathbf{c}=\tfrac12$, one has  
 \begin{align*}%
 	\hat{a}_{42}=\tfrac1{2c_2}-c_3\quad\text{and}\quad 
 	\brat{1-c_2}c_2+\brat{\hat{a}_{42}+c_2 }c_3+c_2=1\,.
 	\end{align*}
 They reduce into the quadratic equation, $c_3^2-\brab{\tfrac1{2c_2}+c_2 }c_3+\brat{c_2-1}^2=0$.
 	If $0<c_2\le\frac{2+\sqrt{6}}{2}$ such that $2 c_2^2- 4 c_2-1\le0$, one gets
 \begin{align}\label{cond: 4-stage IERK2-R}
 c_3^*:=\tfrac{3}{2}\quad\text{if $c_2=1$}\quad\text{and}\quad	c_3^*:=\tfrac1{4c_2}+\tfrac{c_2}{2}-\tfrac{\sqrt{(-2 c_2^2+4 c_2+1) (6 c_2^2-4c_2+1)}}{4c_2}>0\quad\text{if $c_2\neq1$,}
 \end{align} 
 where the other positive root $c_3=\tfrac1{4c_2}+\tfrac{c_2}{2}+\tfrac{\sqrt{(-2 c_2^2+4 c_2+1) (6 c_2^2-4c_2+1)}}{4c_2}$ is dropped for $c_2\neq1$.
 
 By using the definition in \eqref{def: DOC dij}, one has
 \begin{align*}
 	D_{\mathrm{R}}^{(2,4)}(c_2)
 	=\begin{pmatrix}
 		\frac{1}{c_{2}} & 0 & 0 \\[2pt]
 		\frac{2 c_{2}-c_{3}^*}{c_{2}^2} & \frac{1}{c_{2}} & 0 \\[4pt]
 		\frac{4 c_{2}^3-2 c_{2}^2-c_2+2 c_{2}^2c_{3}^*+c_{3}^* -2c_{2} (c_{3}^*)^2}{2 c_{2}^4} & \frac{2c_{2}^2+2c_{3}^* c_{2}-1}{2c_{2}^3} & \frac{1}{c_{2}}
 	\end{pmatrix}.
 \end{align*}
 Now consider the positive definiteness of $D_{\mathrm{R}}^{(2,4)}(c_2)$.
 It is easy to check that the first leading principal minor is positive, the second principal minor $\mathrm{Det}\kbrab{\mathcal{S}(D_{\mathrm{R},2}^{(2,4)};c_2)}=\tfrac{c_{3} (4 c_{2}-c_{3})}{4 c_{2}^4}>0$ for $0<c_2\le\tfrac{2+\sqrt{6}}{2}$,
 while  
 $\mathrm{Det}\kbrab{\mathcal{S}(D_{\mathrm{R}}^{(2,4)};c_2)}=\tfrac{1}{8 c_{2}^7}(c_{2}-1)^2 (2 c_{2}^2+2 c_{2}-1)>0$ for $c_2\ge\tfrac{\sqrt{3}-1}{2}$ and $c_2\neq1$.
 Thus the differentiation matrix $D_{\mathrm{R}}^{(2,4)}(c_2)$ is positive definite if $\tfrac{\sqrt{3}-1}{2}\le c_2<1$ or $1<c_2\le\tfrac{2+\sqrt{6}}{2}$.   
 We obtain the one-parameter R-IERK(2,4;$c_2$) methods with the associated Butcher tableaux
 \begin{equation}\label{scheme: R-IERK-2-4-c2}
  \begin{array}{c|cccc}
 		0 & 0 &  &    &  \\[1pt]
 		c_{2} & \tfrac{c_{2}}{2} & \tfrac{c_{2}}{2} &   \\[1pt]
 		c_3^* & \tfrac{c_3^*-c_{2}}{2} & \tfrac{c_3^*}{2} & \tfrac{c_{2}}{2}    \\[1pt]
 		1 & \frac{1+c_3^*-c_{2}}{2}-\tfrac1{4c_2} & \tfrac{1-c_{2} }{2} & \tfrac1{4c_2}+\tfrac{c_{2}-c_3^*}{2} &   \tfrac{c_{2}}{2} \\[1pt]
 		\hline\\[-10pt]   & \frac{1+c_3^*-c_{2}}{2}-\tfrac1{4c_2} & \tfrac{1-c_{2} }{2} & \tfrac1{4c_2}+\tfrac{c_{2}-c_3^*}{2} &   \tfrac{c_{2}}{2}
 	\end{array},\quad
 	\begin{array}{c|cccc}
 		0 & 0 &  &   &  \\[1pt]
 		c_{2} & c_2 & 0 & &  \\[1pt]
 		c_3^* & c_3^*-c_2 &c_2 & 0&   \\[1pt]
 		1 & 1-\tfrac1{2c_2}+c_3^*-c_2  &\tfrac1{2c_2}-c_3^* & c_2 & 0 \\[1pt]
 		\hline\\[-10pt]   & 1-\tfrac1{2c_2}+c_3^*-c_2 &\tfrac1{2c_2}-c_3^* & c_2&0
 	\end{array}\;,
 \end{equation}
 where $c_3^*$ is defined by \eqref{cond: 4-stage IERK2-R}.
 The associated average dissipation rate is
 \begin{align}\label{rate: IERK2-R}
 	\mathcal{R}_{\mathrm{R}}^{(2,4)}(c_2)=\mathrm{tr}\kbrab{D_{\mathrm{R}}^{(2,4)}(c_2)}=\tfrac{1}{c_2}\quad\text{for $\tfrac{\sqrt{3}-1}{2}\le c_2<1$ or $1<c_2\le\tfrac{2+\sqrt{6}}{2}$}.  
 \end{align} 
 
 \begin{remark}\label{remark: zero-eigenvalue c2_1}
 The case $c_2=1$ in \eqref{rate: IERK2-R} will arrive at $\mathcal{R}_{\mathrm{R}}^{(2,4)}(1)=1$, the optimal value of average dissipation rate. It is noticed that the symmetric matrix $\mathcal{S}(D_{\mathrm{R}}^{(2,4)};1)$ is singular since $\mathrm{Det}\kbrab{\mathcal{S}(D_{\mathrm{R}}^{(2,4)};1)}=0$. 
 This case is permitted in Lemma \ref{lemma: IERK energy stability}, which is valid under the global Lipschitz continuity assumption of nonlinear bulk, cf. \cite[Theorem 2.1 and Corollary 2.1]{LiaoWangWen:2024arxiv}; however, this case should be excluded in our current discussions because it does not fulfill the conditions of Lemma \ref{lemma: convolution terms} and Theorems \ref{thm: regularity of time-discrete solution}-\ref{thm: RIERK convergence}, which require that the symmetric part $\mathcal{S}(D_{\mathrm{R}}^{(2,4)}(c_2))$ has a finite minimum eigenvalue $\lambda_{\min}>0$. It is to mention that, the corresponding R-IERK(2,4;1) scheme still performs well in our numerical experiments, see Section 5, although we cannot theoretically verify the stability and convergence yet.
 \end{remark}

\subsection{Third-order R-IERK methods}

It is not difficult to check that no 4-stage R-IERK 
method exists for the third-order accuracy. One can consider the five-stage R-IERK method
which has only nine independent coefficients.
From Table \ref{table: order condition}, eight order conditions are required such that there remains one independent variable. Unfortunately, we are not able to find any third-order R-IERK methods within five stages to ensure the positive (semi-)definite of the corresponding differentiation matrix  $D_{\mathrm{R}}=D_{\mathrm{E}}$.
We consider six-stage method with one free variable $\hat{a}_{52}$ by choosing $c_{2}=1$, $c_{3}=\tfrac4{5}$,  $c_4=\frac{7}{10}$, $c_5=\frac{12}{25}$ and $\hat{a}_{32}=\frac{6}5$ after a trial and error process. It arrives at  one-parameter R-IERK(3,6; $\hat{a}_{52}$) methods with the following Butcher tableaux
{\small\begin{subequations}\label{scheme: R-IERK-3-6-ah52}
		\begin{align}
			\begin{array}{c|c}
				\mathbf{c} & A \\
				\hline\\[-10pt] 	& \mathbf{b}^T
			\end{array}=
			&\,	\begin{array}{c|cccccc}
				0           & 0 &  &    &  &\\[2pt]
				1           & \tfrac{1}{2} & \tfrac{1}{2} &   &\\[4pt]
				\tfrac{4}{5} & 
				-\tfrac{1}{5} & \tfrac{2}{5} & \tfrac{3}{5}&  \\[4pt]
				\tfrac{7}{10} &					-\tfrac{98716201}{329771428} & -\tfrac{2367068609}{14839714260} & \tfrac{267670251}{412214285} & \tfrac{378048430}{741985713} \\[4pt]
				\tfrac{12}{25} &\frac1{2}\hat{a}_{51}
				& \frac{\hat{a}_{51}+\hat{a}_{52}}{2} & \frac{\hat{a}_{52}+\hat{a}_{53}}{2} & \frac{\hat{a}_{53}+\hat{a}_{54}}{2} & \tfrac{35 \hat{a}_{52}}{103}+\tfrac{89552314349}{363926325900}\\[4pt]
				1            & 	\tfrac{206221}{2306412} & \tfrac{355973}{1153206} & -\tfrac{422005}{768804} & \tfrac{13730}{192201} & \tfrac{553250}{576603} & \tfrac{69125}{576603} \\[4pt]
				\hline\\[-10pt]   & \tfrac{206221}{2306412} & \tfrac{355973}{1153206} & -\tfrac{422005}{768804} & \tfrac{13730}{192201} & \tfrac{553250}{576603} & \tfrac{69125}{576603}
			\end{array},\label{scheme: R-IERK-3-6-ah52 implicit}\\
			\begin{array}{c|c}
				\hat{\mathbf{c}} & \widehat{A} \\
				\hline\\[-10pt] 	& \hat{\mathbf{b}}^T
			\end{array}	=
			&\,\begin{array}{c|cccccc}
				0             & 0 &  &   & & \\[2pt]
				1             & 1 & 0 & & & \\[4pt]
				\tfrac{4}{5} & 	-\tfrac{2}{5} & \tfrac{6}{5} & 0&  & \\[4pt]
				\tfrac{7}{10} &	-\tfrac{98716201}{164885714} & \tfrac{1037580218}{3709928565} & \tfrac{756096860}{741985713}&0  & \\[4pt]
				\tfrac{12}{25} &\hat{a}_{51} & \hat{a}_{52} & \tfrac{14091138538}{30327193825}-\tfrac{190 \hat{a}_{52}}{103} & \tfrac{70 \hat{a}_{52}}{103}+\tfrac{89552314349}{181963162950} &0\\[4pt]
				1            &\tfrac{206221}{1153206} & \tfrac{168575}{384402} & -\tfrac{98430}{64067} & \tfrac{322750}{192201} & \tfrac{138250}{576603}&0\\[4pt]
				\hline\\[-10pt]   & \tfrac{206221}{1153206} & \tfrac{168575}{384402} & -\tfrac{98430}{64067} & \tfrac{322750}{192201} & \tfrac{138250}{576603}&0
			\end{array},\label{scheme: R-IERK-3-6-ah52 explicit}
		\end{align}
\end{subequations}}
where the coefficient $\hat{a}_{51}=\tfrac{17\hat{a}_{52}}{103}-\tfrac{86756827361}{181963162950}$.
Simple calculations yield the differentiation  matrix $D_{\mathrm{R}}^{(3,6)}(\hat{a}_{52})=A_{\mathrm{E}}^{-1}(\hat{a}_{52})E_{5}$, 
which is always
positive definite if $0.664767<\hat{a}_{52}<0.751947$.
The associated average dissipation rate reads
\begin{align}\label{rate: R-IERK-3-6-ah52}
	\mathcal{R}_{\mathrm{R}}^{(3,6)}(\hat{a}_{52})
	=&\,\tfrac{36392632590}{123664285500\hat{a}_{52}+89552314349}
	+\tfrac{219055887768899}{156795586342500}
	\approx \tfrac{1.39708\hat{a}_{52}+1.30599}{\hat{a}_{52}+0.724157}.
\end{align}  
  
We are not sure whether there exists a six-stage IERK scheme having the optimal value of average dissipation rate, $\mathcal{R}_{\mathrm{R}}^{(3,6)}=1$.
Here we do not attempt to construct a fourth-order R-IERK method. Actually, the large number (twenty-eight) of order conditions, cf. \cite{AscherRuuthSpiteri:1997,IzzoJackiewicz:2017,LiaoWangWen:2024arxiv}, makes the fourth-order R-IERK method require at least nine stages so that it would be computationally intensive.

\section{Energy dissipation laws and $L^2$ norm convergence}
\label{sec: convergence}
\setcounter{equation}{0}

\subsection{Technical lemmas and time-discrete system}

We use  the standard seminorms and norms in the Sobolev space $H^{m}(\Omega)$ for $m\ge0$.
Let $C^{\infty}(\Omega)$ be a set of infinitely differentiable $L$-periodic functions defined on $\Omega$,
and $H_{per}^{m}(\Omega)$ be the closure of $C^{\infty}(\Omega)$ in $H^{m}(\Omega)$,
endowed with the semi-norm $|\cdot|_{H_{per}^m}$ and the norm $\mynorm{\cdot}_{H_{per}^{m}}$. 
For simplicity, denote $|\cdot|_{H^m}:=|\cdot|_{H_{per}^m}$, $\mynorm{\cdot}_{H^{m}}:=\mynorm{\cdot}_{H_{per}^{m}}$, $\mynorm{\cdot}_{L^{2}}:=\mynorm{\cdot}_{H^{0}}$, 
and the maximum norm $\mynorm{\cdot}_{L^{\infty}}$.

Next lemma lists some approximations of the $L^2$-projection operator $P_{M}$ and
trigonometric interpolation operator $I_{M}$ defined in subsection 2.1.
\begin{lemma}\cite{ShenTangWang:2011Spectral}\label{lem:Projection-Estimate}
	For any $w\in{H_{per}^{q}}(\Omega)$ and $0\le{\ell}\le{q}$, it holds that
	\begin{align*}
		\mynorm{P_{M}w-w}_{H^{\ell}}
		\le \Ck_wh^{q-s}|w|_{H^{q}},
		\quad \mynorm{P_{M}w}_{H^{\ell}}\le \Ck_w\mynorm{w}_{H^{\ell}};
	\end{align*}
	and, in addition if $q>3/2$,
	\begin{align*}
		\mynorm{I_{M}w-w}_{H^{\ell}}
		\le \Ck_wh^{q-\ell}|w|_{H^{q}},
		\quad \mynorm{I_{M}w}_{H^{\ell}}\le \Ck_w\mynorm{w}_{H^{\ell}}.
	\end{align*}
\end{lemma}

In deriving the stability and error estimate of the stage solutions for semi-implicit multi-stage methods, we need the following  Gr\"{o}nwall-type lemma.
\begin{lemma}\label{lemma: gronwall multi-stage RK1}
	For any integers $s>1$ and $N>1$, let $v_{n,i}$, $g_{n,i}$ and $\tilde{g}_{n,i}$ be the nonnegative discrete functions defined at the stage $t_{n,i}$ for the time level indexes $n=1,2,\cdots,N$  and the stage indexes $i=1,2,\cdots,s$ with a finite $T=t_{N}=t_{N,s}$. Also, let $v_n:=v_{n,s}=v_{n+1,1}$. Assume further that $g_{n,i}$ and $\tilde{g}_{n,i}$ are bounded for $1\le i\le s$ and $n\ge1$.  There exists a positive constant $\omega_0$, independent of the time-step sizes $\tau_n$, such that
	\begin{equation}\label{lemProof: assumption}
		v_{n,i}\le v_{n,1}+\omega_0\tau_n
		\sum_{j=1}^{i-1}v_{n,j}
		+\sum_{j=2}^{i}g_{n,j}+\tilde{g}_{n,i}\quad\text{for $n\ge1$ and $2\le i\le s$.}
	\end{equation} 
	If the maximum step-size $\tau$ is small such that $\tau\le 1/\omega_0$, it holds  that
	\begin{align*}
		v_{n,i}\le&\, 4^{s-1}
		\exp\brab{2^{s-1}\omega_0t_{n-1}}
		\braB{v_0+\sum_{k=1}^{n-1}\sum_{j=2}^{s}g_{k,j}+\sum_{j=2}^{i}g_{n,j}
			+\sum_{k=1}^{n}\max_{2\le j\le s}\tilde{g}_{k,j}}
	\end{align*}
	for $1\le n\le N$ and $2\le i\le s$.
\end{lemma}
\begin{proof}For fixed $n$, a simple induction for \eqref{lemProof: assumption} gives the stage estimate
	\begin{align}\label{lemProof: internal error}
		v_{n,i}\le&\, \bra{1+\omega_0\tau_n}^{i-1}
		\braB{v_{n,1}+\sum_{\ell=2}^{i}g_{n,\ell}+\max_{2\le j\le s}\tilde{g}_{n,j}}		
	\end{align}
	for $2\le i\le s$. Taking 	$i=s$ and $n=k$ in the inequality \eqref{lemProof: assumption} yields
	\begin{equation}
		v_{k}=v_{k,s}\le v_{k-1}+\omega_0\tau_k
		\sum_{j=1}^{s-1}v_{k,j}
		+\sum_{j=2}^{s}g_{k,j}+\tilde{g}_{k,s}\quad\text{for $k\ge1$,}
	\end{equation} 
	and, by summing the index $k$ from $1$ to $n$,
	\begin{equation*}
		v_{n}\le v_0+\omega_0	\sum_{k=1}^{n}
		\sum_{j=1}^{s-1}\tau_kv_{k,j}
		+\sum_{k=1}^{n}\sum_{j=2}^{s}g_{k,j}+\sum_{k=1}^{n}\tilde{g}_{k,s}
		\quad\text{for $n\ge1$.}
	\end{equation*}
	Now we insert the stage estimate \eqref{lemProof: internal error} into the above inequality to obtain that
	\begin{align}\label{lemProof: level error}
		v_{n}
		\le&\, \bra{1+\omega_0\tau_1}^{s-1}v_0+\omega_0	\sum_{k=2}^{n}\tau_kv_{k-1}
		\sum_{j=1}^{s-1}
		\bra{1+\omega_0\tau_k}^{j-1}	
		\nonumber\\&\,
		+	\sum_{k=1}^{n}\bra{1+\omega_0\tau_k}^{s-1}\sum_{j=2}^{s}g_{k,j}	
		+\sum_{k=1}^{n}\bra{1+\omega_0\tau_k}^{s-1}\max_{2\le j\le s}\tilde{g}_{k,j}
	\end{align}
	for $n\ge1$, where we used the facts $g_{k,j}\ge0$, $\tilde{g}_{k,j}\ge0$ and $1+\omega_0\tau_k\sum_{j=1}^{s-1}\bra{1+\omega_0\tau_k}^{j-1}
	=\bra{1+\omega_0\tau_k}^{s-1}.$
	Under the maximum step-size setting $\tau\le 1/\omega_0$, one has
	$$\bra{1+\omega_0\tau_k}^{s-1}\le 2^{s-1}\quad\text{and}\quad\sum_{j=1}^{s-1}\bra{1+\omega_0\tau_k}^{j-1}
	\le \sum_{j=1}^{s-1}2^{j-1}= 2^{s-1}-1.$$
	It follows from \eqref{lemProof: level error} that
	\begin{align*}
		v_{n}\le&\, 2^{s-1}v_0+2^{s-1}\omega_0	\sum_{k=2}^{n}\tau_kv_{k-1}
		+2^{s-1}	\sum_{k=1}^{n}\sum_{j=2}^{s}g_{k,j}+2^{s-1}\sum_{k=1}^{n}\max_{2\le j\le s}\tilde{g}_{k,j}\quad\text{for $n\ge1$.}
	\end{align*}
	The discrete Gr\"{o}nwall inequality, cf. \cite[Lemma 3.1]{LiaoZhang:2021}, leads to
	the time-level estimate,
	\begin{align*}
		v_{n}\le&\, 2^{s-1}
		\exp\brab{2^{s-1}\omega_0t_{n}}
		\braB{v_0+\sum_{k=1}^{n}\sum_{j=2}^{s}g_{k,j}+\sum_{k=1}^{n}\max_{2\le j\le s}\tilde{g}_{k,j}}\quad\text{for $1\le n\le N$.}
	\end{align*}
	Return to the stage estimate \eqref{lemProof: internal error} with the setting $\tau\le 1/\omega_0$, we get the time-stage estimate
	\begin{align*}
		v_{n,i}\le&\, 4^{s-1}
		\exp\brab{2^{s-1}\omega_0t_{n-1}}
		\braB{v_0+\sum_{k=1}^{n-1}\sum_{j=2}^{s}g_{k,j}+\sum_{j=2}^{i}g_{n,j}
			+\sum_{k=1}^{n-1}\max_{2\le j\le s}\tilde{g}_{k,j}+\max_{2\le i\le s}\tilde{g}_{n,j}}
	\end{align*}
	for $1\le n\le N$ and $1\le i\le s$. It completes the proof.
\end{proof}

To establish the original energy dissipation laws and the  $L^2$ norm error estimate of the fully discrete R-IERK method \eqref{Scheme: R-IERK stabilized} for the fourth-order nonlinear CH model, we will update the idea of time-space error splitting approach in \cite{LiGaoSun:2014,LiSun:2013} via the following three steps: 
\begin{itemize}[itemindent=5mm]
	\item[(Step 1)] Subsection 4.2 addresses the time-discrete system to establish the regularity of time-discrete stage solutions $U^{n,i}$ via the energy arguments with a rough setting of stage defects. 
	Always, we assume that the initial data and the solution of \eqref{cont: Problem-CH} are smooth. There exist two integers $m\ge1$, $p\ge1$ and a constant $\ck_{\phi}>0$ such that
	\begin{align}\label{regular: exact solution}				
		\mynormb{\Phi^0}_{H^{m+4}}
		+\sum_{k=0}^{2}\mynormb{\partial_t^{(k)}\Phi(t)}_{H^{m+4-k}}
		+\sum_{k=3}^{p+1}\mynormb{\partial_t^{(k)}\Phi(t)}_{L^{2}}\le \ck_{\phi}\quad\text{for $0<t<T$}.
	\end{align}
	\item[(Step 2)] Subsection 4.3 addresses the fully-discrete $L^2$ norm error estimate by using the pseudo-spectral approximation of time-discrete system and establishes the maximum norm boundedness of stage solutions $u_h^{n,i}$ via the inverse estimate such that  we can establish the original energy dissipation law of the fully discrete R-IERK method \eqref{Scheme: R-IERK stabilized}.
	\item[(Step 3)] With the established regularity of time-discrete solution and the maximum norm boundedness of stage solutions, Subsection 4.4 arrives at the $L^2$ norm error estimate for the fully discrete R-IERK method \eqref{Scheme: R-IERK stabilized} with the full accuracy.
\end{itemize}

To process the numerical analysis of the fully discrete R-IERK method \eqref{Scheme: R-IERK stabilized}, we will consider time-discrete stage solution $U^{n,i}$ satisfying the following time-discrete system 
\begin{align}\label{system: time discrete R-IERK}				
	U^{n,i+1}=&\,U^{n,1}+\tau_n\sum_{j=1}^{i}\hat{a}_{i+1,j}\Delta \kbraB{L_{\kappa}U^{n,j+\frac1{2}}-f_{\kappa}(U^{n,j})}
	\quad\text{for $n\ge1$ and $1\le i\le s_{\mathrm{I}}$,}
\end{align}
subject to the initial data $U^{1,1}=\Phi^{1,1}:=\Phi(\myvec{x},0)$.
Inserting the exact solution at the stage $t_{n,i}$, $\Phi^{n,i}:=\Phi(\myvec{x},t_{n-1}+c_{i}\tau_n)$, into the time-discrete system, one has
\begin{align}\label{system: time approximation R-IERK}				
	\Phi^{n,i+1}=&\,\Phi^{n,1}+\tau_n\sum_{j=1}^{i}\hat{a}_{i+1,j}\Delta \kbraB{L_{\kappa}\Phi^{n,j+\frac1{2}}-f_{\kappa}(\Phi^{n,j})}
	+\tau_n\zeta_{\mathrm{R}}^{n,i+1}
	\quad\text{for $n\ge1$ and $1\le i\le s_{\mathrm{I}}$,}
\end{align}
where the temporal defects $\zeta_{\mathrm{R}}^{n,i+1}$ at the stage $t_{n,i+1}$ will be determined later (with $\zeta_{\mathrm{R}}^{n,1}\equiv0$). The time errors $\tilde{U}^{n,i+1}:=\Phi^{n,i+1}-U^{n,i+1}$ fulfill the time-error system
\begin{align}\label{system: time error R-IERK}				
	\tilde{U}^{n,i+1}=&\,\tilde{U}^{n,1}+\tau_n\sum_{j=1}^{i}\hat{a}_{i+1,j}\Delta \kbraB{L_{\kappa}\tilde{U}^{n,j+\frac1{2}}
		+f_{\kappa}(U^{n,j})-f_{\kappa}(\Phi^{n,j})}
	+\tau_n\zeta_{\mathrm{R}}^{n,i+1},
\end{align}
or, with the matrix $(\underline{\hat{a}}_{i+1,j})_{s_{\mathrm{I}}\times s_{\mathrm{I}}}=E_{s_{\mathrm{I}}}^{-1}A_{\mathrm{E}}=D_{\mathrm{R}}^{-1}$,
\begin{align}\label{system: time error R-IERK-difference}				
	\delta_{\tau}\tilde{U}^{n,i+1}=\tau_n\sum_{j=1}^{i}\underline{\hat{a}}_{i+1,j}\Delta \kbraB{L_{\kappa}\tilde{U}^{n,j+\frac1{2}}
		+f_{\kappa}(U^{n,j})-f_{\kappa}(\Phi^{n,j})}
	+\tau_n\delta_{\tau}\zeta_{\mathrm{R}}^{n,i+1}
\end{align} 
for $n\ge1$ and $1\le i\le s_{\mathrm{I}}$.  Note that, the exact solution $\Phi^{n,i}$ and the time-discrete solution $U^{n,i}$ preserve the volume conservation such that $\brab{\tilde{U}^{n,i},1}=0$ and $\brab{\zeta_{\mathrm{R}}^{n,i},1}=0$ for $n\ge1$ and $1\le i\le s$.
Recalling the  orthogonal property, cf. \cite[subsection 2.1]{LiaoWangWen:2024arxiv}, one has
$$\sum_{i=1}^{k}	d^{(R)}_{k,i}\sum_{j=1}^{i}\underline{\hat{a}}_{i+1,j}v^j
=\sum_{j=1}^{k}v^j\sum_{i=j}^{k}d^{(R)}_{k,i}\underline{\hat{a}}_{i+1,j}\equiv v^k
\quad \text{for $1\le k\le s_{\mathrm{I}}$.}$$
Multiplying both sides of the equality \eqref{system: time error R-IERK-difference}	 by the kernels $d^{(R)}_{j,i}$ defined in \eqref{def: DOC dij}, and summing the stage index $i$ from $i=1$ to $j$, it is easy to obtain the following equivalent form of \eqref{system: time error R-IERK-difference},
\begin{align}\label{system: DOC time error R-IERK-difference}					
	\sum_{i=1}^{j}d^{(R)}_{j,i}\delta_{\tau}\tilde{U}^{n,i+1}=&\,\tau_n\Delta \kbraB{L_{\kappa}\tilde{U}^{n,j+\frac1{2}}
		+f_{\kappa}(U^{n,j})-f_{\kappa}(\Phi^{n,j})}
	+\tau_n\sum_{i=1}^{j}d^{(R)}_{j,i}\delta_{\tau}\zeta_{\mathrm{R}}^{n,i+1}
\end{align} 
for $n\ge1$ and $1\le j\le s_{\mathrm{I}}$. 

It is to mention that, the time-error system \eqref{system: time error R-IERK-difference}	will be used for the $L^2$ norm estimate of stage solution; while the equivalent form \eqref{system: DOC time error R-IERK-difference} will be used for the $H_{per}^{2}$ seminorm estimate. To treat with the involved method coefficients $\underline{\hat{a}}_{i+1,j}$ and $d^{(R)}_{i,j}$, we will use the following lemma.
\begin{lemma}\label{lemma: convolution terms}	
	Assume that the differentiation matrix $D_{\mathrm{R}}=A_{\mathrm{E}}^{-1}E_{s_{\mathrm{I}}}$ of the R-IERK method \eqref{Scheme: R-IERK stabilized} and the inverse $D_{\mathrm{R}}^{-1}$ are positive definite.
	Let $\lambda_{\min}$ and $\sigma_{\min}$ be the minimum eigenvalues of the symmetric parts $\mathcal{S}(D_{\mathrm{R}})$ and $\mathcal{S}(D_{\mathrm{R}}^{-1})$, respectively. Then, for any time sequences $\{v^j,u^j\,|\,j\ge1\}$, it holds that
	\begin{itemize}
		\item[(i)] $\sum\limits_{i=1}^{k}\sum\limits_{j=1}^{i}d^{(R)}_{i,j}v^jv^i\ge 
		\lambda_{\min}\sum\limits_{i=1}^{k}\brat{v^i}^2$,\;\; and \;\; $\sum\limits_{i=1}^{k}\sum\limits_{j=1}^{i}
		\underline{\hat{a}}_{i+1,j}v^jv^i\ge 
		\sigma_{\min}\sum\limits_{i=1}^{k}\brat{v^i}^2$;
		\item[(ii)] $\sum\limits_{i=1}^{k}\sum\limits_{j=1}^{i}d^{(R)}_{i,j}v^ju^i\le 
		\frac1{\sigma_{\min}}\sum\limits_{i=1}^{k}\absb{v^i}\abs{u^i}$,\;\; and \;\; $\sum\limits_{i=1}^{k}\sum\limits_{j=1}^{i}\underline{\hat{a}}_{i+1,j}v^ju^i\le 
		\frac1{\lambda_{\min}}\sum\limits_{i=1}^{k}\absb{v^i}\absb{u^i}$.
	\end{itemize}
\end{lemma}
\begin{proof}The result (i) is obvious. Since $D_{\mathrm{R}}-\lambda_{\min}I$ is positive semi-definite, \cite[Lemma 5.2]{ChenYuZhangZhang:2024} says that the maximum eigenvalue of $(D_{\mathrm{R}}^{-1})^TD_{\mathrm{R}}^{-1}$ can be bounded by $1/\lambda_{\min}^2$. Similarly, the maximum eigenvalue of $D_{\mathrm{R}}^TD_{\mathrm{R}}$  can be bounded by $1/{\sigma_{\min}^2}$. They imply the result (ii) and complete the proof.
\end{proof}

In handling the nonlinear term, the following technical lemma will be helpful, while the proof is simple due to the formula 
$F'(v)-F'(w)=(v-w)\int_0^1F''[\gamma v+(1-\gamma) w]\zd\gamma$.
Actually, it is easy to derive that
$F'(\Phi^{n,j})-F'(U^{n,j})=\tilde{U}^{n,j}\int_0^1F''\brab{U^{n,j}+\gamma \tilde{U}^{n,j}}\zd\gamma$,
and 
\begin{align*}
	&\,\delta_{\tau}F'(\Phi^{n,j})-\delta_{\tau}F'(U^{n,j})
	=\delta_{\tau}\tilde{U}^{n,j}
	\int_0^1F''\brab{U^{n,j-1}+\beta \delta_{\tau}U^{n,j}}\zd\beta\\
	&\,\hspace{1.5cm}+\delta_{\tau}\Phi^{n,j}\int_0^1\brab{\tilde{U}^{n,j-1}
		+\beta \delta_{\tau}\tilde{U}^{n,j}}\int_0^1F'''\kbraB{U^{n,j-1}
		+\beta \delta_{\tau}U^{n,j}+\gamma\brab{\tilde{U}^{n,j-1}
			+\beta \delta_{\tau}\tilde{U}^{n,j}}}\zd\gamma\zd\beta,
\end{align*}
where $\delta_{\tau}\Phi^{n,j}=(t_{n,j}-t_{n,j-1})\int_0^1\frac{\zd}{\zd t}\Phi[t_{n,j-1}+\beta(t_{n,j}-t_{n,j-1})]\zd\beta$.
\begin{lemma}\label{lemma: nonlinear term}
	Under the regularity setting \eqref{regular: exact solution}, it holds that
	\begin{itemize}
		\item[(i)] $\mynormb{F'(\Phi^{n,j})-F'(U^{n,j})}_{L^2}\le\mynormb{\tilde{U}^{n,j}}_{L^2}\int_0^1\mynormb{F''\brab{\mathrm{Q}_{\gamma}^{n,j}[U]}}_{L^\infty}\zd\gamma$, where
		$\mathrm{Q}_{\gamma}^{n,j}[U]:=U^{n,j}+\gamma \tilde{U}^{n,j}$;
		\item[(ii)] $\mynormb{\delta_{\tau}\Phi^{n,j}}_{L^\infty}\leq \abs{t_{n,j}-t_{n,j-1}}\mynorm{\partial_t\Phi(t)}_{L^\infty}
		\le \ck_{\phi}\tau_n$ and then
		\begin{align*}
			&\,\mynormb{\delta_{\tau}F'(\Phi^{n,j})-\delta_{\tau}F'(U^{n,j})}_{L^2}
			\le\mynormb{\delta_{\tau}\tilde{U}^{n,j}}_{L^2}
			\int_0^1\mynormb{F''\brab{\mathsf{Q}_{\beta}^{n,j}[U]}}_{L^\infty}\zd\beta\\
			&\,\hspace{2cm}+2\ck_{\phi}\tau_n\brab{\mynormb{\tilde{U}^{n,j}}_{L^2}
				+\mynormb{\tilde{U}^{n,j-1}}_{L^2}}
			\int_0^1\int_0^1\mynormb{F'''\brab{\mathcal{Q}_{\beta,\gamma}^{n,j}[U]}}_{L^\infty}
			\zd\gamma\zd\beta,
		\end{align*}
		where $\mathsf{Q}_{\beta}^{n,j}[U]:=U^{n,j-1}+\beta \delta_{\tau}U^{n,j}$ and
		$\mathcal{Q}_{\beta,\gamma}^{n,j}[U]:=U^{n,j-1}
		+\beta \delta_{\tau}U^{n,j}+\gamma\brab{\tilde{U}^{n,j-1}
			+\beta \delta_{\tau}\tilde{U}^{n,j}}$.
	\end{itemize}	
\end{lemma}

\subsection{Stage regularity of time-discrete solutions}

 To establish the regularity of time stage solutions $U^{n,i}$ to the time-discrete system \eqref{system: time discrete R-IERK}, we impose the following rough assumption
\begin{align}\label{system: rough assumption}
	\mynormb{\zeta_{\mathrm{R}}^{n,i+1}}_{L^2}\le \ck_1\tau_n
	\quad \text{for $1\le n\le N$ and $1\le i\le s_{\mathrm{I}}$.}
\end{align}  
This setting \eqref{system: rough assumption} is imposed according to two facts: (1) assuming $\zeta_{\mathrm{R}}^{n,i+1}=0$ for $1\le i\le s_{\mathrm{I}}-1$ would be not reasonable to derive the errors of stage solutions; (2) the first-order setting is to make the present analysis extendable to the Radau-type IERK methods in \cite{FuTangYang:2024,LiaoWangWen:2024arxiv}.

\subsubsection{$H_{per}^{2}$ norm boundedness of stage solutions}
We will use the complete mathematical induction to prove that the stage solutions $U^{n,i}$ of the time-discrete system \eqref{system: time discrete R-IERK} are bounded in the $H_{per}^2$ norm, that is,
\begin{align}\label{bound: H2 time error}
	\mynormb{\tilde{U}^{n,i}}_{H^2}\le \ck_{10}/\epsilon^2\quad\text{for $1\le n\le N$ and $1\le i\le s$.}
\end{align}

Obviously, it holds for $n=1$ and $i=1$ since $U^{1,1}=\Phi^0$. Put the inductive hypothesis
\begin{align}\label{bound: H2 time error hypothesis}
	\mynormb{\tilde{U}^{l,\ell}}_{H^2}\le \ck_{10}/\epsilon^2\quad\text{such that}\quad
	\mynormb{U^{l,\ell}}_{L^\infty}\le \ck_{10}^*/\epsilon^2
	\quad\text{for $1\le l\le n$ and $1\le \ell\le k$,}
\end{align}
where $\ck_{10}^*=\ck_{\phi}\epsilon^2+\ck_{\Omega}\ck_{10}$.
Now we are to prove that $\mynormb{\tilde{U}^{n,k+1}}_{H^2}\le \ck_{10}/\epsilon^2$.

\paragraph{($L^2$ norm rough estimate)}
Making the $L^2$ inner product of \eqref{system: time error R-IERK-difference} by $2\tilde{U}^{n,i+\frac{1}{2}}$, applying the discrete Green's formula and summing the stage index $i$ from $i=1$ to $k$, one can find that
\begin{align}\label{system: L2 time error R-IERK-difference}
	\mynormb{\tilde{U}^{n,k+1}}_{L^2}^2-\mynormb{\tilde{U}^{n,1}}_{L^2}^2
	=&\,2\tau_n\sum_{i=1}^{k}\sum_{j=1}^{i}\underline{\hat{a}}_{i+1,j}
	\brab{
		L_{\kappa}\tilde{U}^{n,j+\frac{1}{2}},\Delta\tilde{U}^{n,i+\frac{1}{2}}}_{L^2}
	+2\tau_n\sum_{i=1}^{k}
	\brab{\delta_{\tau}\zeta_{\mathrm{R}}^{n,i+1},\tilde{U}^{n,i+\frac{1}{2}}}_{L^2}\notag\\
	&\,+2\tau_n\sum_{i=1}^{k}\sum_{j=1}^{i}
	\underline{\hat{a}}_{i+1,j}\brab{f_{\kappa}(U^{n,j})-f_{\kappa}(\Phi^{n,j}), \Delta\tilde{U}^{n,i+\frac{1}{2}}}_{L^2}.
\end{align}
Lemma \ref{lemma: convolution terms} (i) yields
\begin{align*}
	2\tau_n\sum_{i=1}^{k}\sum_{j=1}^{i}\underline{\hat{a}}_{i+1,j}\brab{
		L_{\kappa}\tilde{U}^{n,j+\frac{1}{2}},\Delta\tilde{U}^{n,i+\frac{1}{2}}}_{L^2}
	\le-2\epsilon^2\sigma_{\min}\tau_n \sum_{i=1}^{k}\mynormb{\Delta\tilde{U}^{n,i+\frac{1}{2}}}_{L^2}^2.
\end{align*}
For the last term in \eqref{system: L2 time error R-IERK-difference}, one can apply Lemma \ref{lemma: convolution terms} (ii) and Lemma \ref{lemma: nonlinear term} (i) to derive that
\begin{align*}
	2\tau_n\sum_{i=1}^{k}\sum_{j=1}^{i}\underline{\hat{a}}_{i+1,j}\brab{
		f_{\kappa}(U^{n,j})-f_{\kappa}(\Phi^{n,j}),\Delta\tilde{U}^{n,i+\frac{1}{2}}}_{L^2}
	\le \frac{2\tau_n}{\lambda_{\min}} \sum_{i=1}^{k}\mynormb{f_{\kappa}(U^{n,i})-f_{\kappa}(\Phi^{n,i})}_{L^2}
	\mynormb{\Delta\tilde{U}^{n,i+\frac{1}{2}}}_{L^2}\\
	\le 
	2\epsilon^2\sigma_{\min}\tau_n
	\sum_{i=1}^{k}\mynormb{\Delta\tilde{U}^{n,i+\frac{1}{2}}}_{L^2}^2
	+\frac{\tau_n}{2\lambda_{\min}^2\sigma_{\min}\epsilon^2}	 \sum_{i=1}^{k}\mynormb{\tilde{U}^{n,i}}_{L^2}^2\int_0^1\brab{\kappa+\mynormb{F''\brab{\mathrm{Q}_{\gamma}^{n,i}[U]}}_{L^\infty}}^2\zd\gamma,
\end{align*}
where $\mathrm{Q}_{\gamma}^{n,i}[U]$ is defined in Lemma \ref{lemma: nonlinear term} and then $\mynormb{F''\brab{\mathrm{Q}_{\gamma}^{n,i}[U]}}_{L^\infty}\le \ck_2:=\max\limits_{\mynormt{\xi}_{L^\infty}\le \ck_{10}^*/\epsilon^2}\mynormb{F''(\xi)}_{L^\infty}$ according to the uniform bound \eqref{bound: H2 time error hypothesis}.
Thus, it follows from \eqref{system: L2 time error R-IERK-difference} that
\begin{align}\label{system: L2 time error R-IERK-difference2}
	\mynormb{\tilde{U}^{n,k+1}}_{L^2}^2-\mynormb{\tilde{U}^{n,1}}_{L^2}^2
	\le\frac{\ck_3\tau_n}{\epsilon^2} \sum_{i=1}^{k}\mynormb{\tilde{U}^{n,i}}_{L^2}^2
	+2\tau_n\sum_{i=1}^{k}\mynormb{\delta_{\tau}\zeta_{\mathrm{R}}^{n,i+1}}_{L^2}\mynormb{\tilde{U}^{n,i+1}}_{L^2},
\end{align}
where the constant $\ck_3:=\frac{(\ck_2+\kappa)^2}{2\sigma_{\min}\lambda_{\min}^2}$. For any index $n$, choose  a finite $k_0$ satisfying $1\le k_0\le k$ such that $\mynormb{\tilde{U}^{n,k_0+1}}_{L^2}:=\max_{0\le \ell\le k}\mynormb{\tilde{U}^{n,\ell+1}}_{L^2}$. We set $k=k_0$ in \eqref{system: L2 time error R-IERK-difference2} to get
\begin{align*}
	\mynormb{\tilde{U}^{n,k_0+1}}_{L^2}\le\mynormb{\tilde{U}^{n,1}}_{L^2}
	+\frac{\ck_3\tau_n}{\epsilon^2} \sum_{i=1}^{k_0}\mynormb{\tilde{U}^{n,i}}_{L^2}	
	+2\tau_n\sum_{i=1}^{k_0}\mynormb{\delta_{\tau}\zeta_{\mathrm{R}}^{n,i+1}}_{L^2},
\end{align*}
and then, due to $k_0\le k$ and $\mynormb{\tilde{U}^{n,k+1}}_{L^2}\le\mynormb{\tilde{U}^{n,k_0+1}}_{L^2}$,
\begin{align}\label{system: L2 time error R-IERK-difference3}
	\mynormb{\tilde{U}^{n,k+1}}_{L^2}\le\mynormb{\tilde{U}^{n,1}}_{L^2}
	+\frac{\ck_3\tau_n}{\epsilon^2}\sum_{i=1}^{k}\mynormb{\tilde{U}^{n,i}}_{L^2}
	+2\tau_n\sum_{i=1}^{k}\mynormb{\delta_{\tau}\zeta_{\mathrm{R}}^{n,i+1}}_{L^2}.
\end{align}
By using the defect bound \eqref{system: rough assumption}, the discrete Gr\"{o}nwall inequality in Lemma \ref{lemma: gronwall multi-stage RK1} 
together with the maximum time-step condition $\tau\le \epsilon^2/\ck_3$ gives
\begin{align}\label{system: L2 time error R-IERK-difference4}
	\mynormb{\tilde{U}^{n,k+1}}_{L^2}\le&\, 4^{s}
	\exp\brab{2^{s-1}\ck_3t_{n-1}/\epsilon^2}s_{\mathrm{I}}\ck_1t_n\tau\le \ck_4\tau,
\end{align}
where the constant $\ck_4:=4^{s}
\exp\brab{2^{s-1}\ck_3T/\epsilon^2}s_{\mathrm{I}}\ck_1T$.

\paragraph{($H_{per}^2$ seminorm rough estimate)}
Making the $L^2$ inner product of \eqref{system: DOC time error R-IERK-difference} by $\frac{2}{\epsilon^2}\delta_{\tau}\tilde{U}^{n,j+1}$, applying the discrete Green's formula and summing $j$ from $j=1$ to $k$, one can find that
\begin{align}\label{system: H2 time error R-IERK-difference}
	&\mynormb{\Delta\tilde{U}^{n,k+1}}_{L^2}^2-\mynormb{\Delta\tilde{U}^{n,1}}_{L^2}^2
	=-\frac2{\epsilon^2\tau_n}\sum_{j=1}^{k}\sum_{i=1}^{j}	d^{(R)}_{j,i}\brab{\delta_{\tau}\tilde{U}^{n,i+1},\delta_{\tau}\tilde{U}^{n,j+1}}_{L^2}
	-\frac{\kappa}{\epsilon^2}\sum_{j=1}^{k}
	\mynormb{\delta_{\tau}\nabla\tilde{U}_h^{n,j+1}}_{L^2}^2\notag\\
&\hspace{1.4cm}	+\frac2{\epsilon^2}\sum_{j=1}^{k}
	\brab{F'(U^{n,j})-F'(\Phi^{n,j}),\delta_{\tau}\Delta\tilde{U}^{n,j+1}}_{L^2}
	+\frac2{\epsilon^2}\sum_{j=1}^{k}\sum_{i=1}^{j}d^{(R)}_{j,i}
	\brab{\delta_{\tau}\zeta_{\mathrm{R}}^{n,i+1},\delta_{\tau}\tilde{U}^{n,j+1}}_{L^2}.
\end{align}
According to Lemma \ref{lemma: convolution terms} (i), one has  
\begin{align*}
	-\frac2{\epsilon^2\tau_n}\sum_{j=1}^{k}\sum_{i=1}^{j}	d^{(R)}_{j,i}\brab{\delta_{\tau}\tilde{U}^{n,i+1}, \delta_{\tau}\tilde{U}^{n,j+1}}_{L^2}
	\le-\frac{2\lambda_{\min}}{\epsilon^2\tau_n} \sum_{j=1}^{k}\mynormb{\delta_{\tau}\tilde{U}^{n,j+1}}_{L^2}^2.
\end{align*}
With the help of Lemma \ref{lemma: convolution terms} (ii) and the Young inequality, one gets
\begin{align*}
	\frac2{\epsilon^2}\sum_{j=1}^{k}\sum_{i=1}^{j}
	d^{(R)}_{j,i}\brab{\delta_{\tau}\zeta_{\mathrm{R}}^{n,i+1},\delta_{\tau}\tilde{U}^{n,j+1}}_{L^2}
	\le&\,\frac{2}{\sigma_{\min}\epsilon^2}\sum_{j=1}^{k}
	\mynormb{\delta_{\tau}\zeta_{\mathrm{R}}^{n,j+1}}_{L^2}
	\mynormb{\delta_{\tau}\tilde{U}^{n,j+1}}_{L^2}\\
	\le &
	\frac{\lambda_{\min}}{\epsilon^2\tau_n} \sum_{j=1}^{k}\mynormb{\delta_{\tau}\tilde{U}^{n,j+1}}_{L^2}^2
	+\frac{\tau_n}{\lambda_{\min}\sigma_{\min}^{2}\epsilon^2}
	\sum_{j=1}^{k}\mynormb{\delta_{\tau}\zeta_{\mathrm{R}}^{n,j+1}}_{L^2}^2.
\end{align*}
To bound the nonlinear term, we recall the Abel-type formula of summation-by-part, also see \cite{LiaoSun:2010,LiSunZhao:2012},  $\sum_{j=1}^{k}u^j(v^{j+1}-v^{j})
=u^{k}v^{k+1}-\sum_{j=2}^{k}v^{j}(u^j-u^{j-1})-u^1v^{1}$.
One applies Lemma \ref{lemma: nonlinear term} to get
\begin{align*}
	J_{k}:=&\,\frac2{\epsilon^2}\sum_{j=1}^{k}
	\brab{F'(U^{n,j})-F'(\Phi^{n,j}),
		\delta_{\tau}\Delta\tilde{U}^{n,j+1}}_{L^2}
	=\frac2{\epsilon^2}\brab{F'(U^{n,k})-F'(\Phi^{n,k}),
		\Delta\tilde{U}^{n,k+1}}_{L^2}\\
	&+\frac2{\epsilon^2}
	\brab{F'(U^{n,1})-F'(\Phi^{n,1}),
		\Delta\tilde{U}^{n,1}}_{L^2}-\frac2{\epsilon^2}\sum_{j=2}^{k}
	\brab{\delta_{\tau}F'(U^{n,j})-\delta_{\tau}F'(\Phi^{n,j}),
		\Delta\tilde{U}^{n,j}}_{L^2}\triangleq\sum_{\ell=1}^3J_{k,\ell},
\end{align*}
where, by using the discrete functionals 
$\mathrm{Q}_{\gamma}^{n,j}[U]$, $\mathsf{Q}_{\beta}^{n,j}[U]$ 
and $\mathcal{Q}_{\beta,\gamma}^{n,j}[U]$ 
defined in Lemma \ref{lemma: nonlinear term},
\begin{align*} J_{k,1}\le&\,\frac2{\epsilon^2}\mynormb{\tilde{U}^{n,k}}_{L^2}
	\mynormb{\Delta\tilde{U}^{n,k+1}}_{L^2}
	\int_0^1\mynormb{F''\brab{\mathrm{Q}_{\gamma}^{n,k}[U]}}_{L^\infty}\zd\gamma,\\
	J_{k,2}\le&\,\frac2{\epsilon^2}\mynormb{\tilde{U}^{n,1}}_{L^2}
	\mynormb{\Delta\tilde{U}^{n,1}}_{L^2}
	\int_0^1\mynormb{F''\brab{\mathrm{Q}_{\gamma}^{n,1}[U]}}_{L^\infty}\zd\gamma,
\\
	J_{k,3}\le&\,\frac{4\ck_{\phi}\tau_n}{\epsilon^2}\sum_{j=2}^{k}
	\mynormb{\Delta\tilde{U}^{n,j}}_{L^2}
	\brab{\mynormb{\tilde{U}^{n,j}}_{L^2}+\mynormb{\tilde{U}^{n,j-1}}_{L^2}}
	\int_0^1\int_0^1\mynormb{F'''\brab{\mathcal{Q}_{\beta,\gamma}^{n,j}[U]}}_{L^\infty}
	\zd\gamma\zd\beta\\
	&+\frac2{\epsilon^2}\sum_{j=2}^{k}
	\mynormb{\Delta\tilde{U}^{n,j}}_{L^2}\mynormb{\delta_{\tau}\tilde{U}^{n,j}}_{L^2}
	\int_0^1\mynormb{F''\brab{\mathsf{Q}_{\beta}^{n,j}[U]}}_{L^\infty}\zd\beta\\
	\le &\frac{4\ck_{\phi}\tau_n}{\epsilon^2}
	\sum_{j=2}^{k}\mynormb{\Delta\tilde{U}^{n,j}}_{L^2}
	\brab{\mynormb{\tilde{U}^{n,j}}_{L^2}+\mynormb{\tilde{U}^{n,j-1}}_{L^2}}
	\int_0^1\int_0^1\mynormb{F'''\brab{\mathcal{Q}_{\beta,\gamma}^{n,j}[U]}}_{L^\infty}
	\zd\gamma\zd\beta\\
	&+\frac{\tau_n}{\lambda_{\min}\epsilon^2}\sum_{j=2}^{k}
	\mynormb{\Delta\tilde{U}^{n,j}}_{L^2}^2
	\int_0^1\mynormb{F''\brab{\mathsf{Q}_{\beta}^{n,j}[U]}}_{L^\infty}^2\zd\beta
	+\frac{\lambda_{\min}}{\epsilon^2\tau_n} \sum_{j=1}^{k}\mynormb{\delta_{\tau}\tilde{U}^{n,j+1}}_{L^2}^2.
\end{align*}
Applying the maximum norm  bound \eqref{bound: H2 time error hypothesis} of stage solutions, one has $\mynormb{F''\brab{\mathrm{Q}_{\gamma}^{n,i}[U]}}_{L^\infty}\le \ck_2,$
$$\mynormb{F''\brab{\mathsf{Q}_{\beta}^{n,i}[U]}}_{L^\infty}\le \ck_2\quad\text{and}\quad
\mynormb{F'''\brab{\mathcal{Q}_{\beta,\gamma}^{n,i}[U]}}_{L^\infty}\le \ck_5:=\max\limits_{\mynormt{\xi}_{L^\infty}\le 2\ck_{10}^*/\epsilon^2}\mynormb{F'''(\xi)}_{L^\infty}.$$
By collecting the above estimates, 
it follows from \eqref{system: H2 time error R-IERK-difference} that 
\begin{align}\label{system: H2 time error R-IERK-difference2}
	\mynormb{\Delta\tilde{U}^{n,k+1}}_{L^2}^2
	\le&\,\mynormb{\Delta\tilde{U}^{n,1}}_{L^2}^2
	+\frac{2\ck_2}{\epsilon^2}\mynormb{\tilde{U}^{n,k}}_{L^2}
	\mynormb{\Delta\tilde{U}^{n,k+1}}_{L^2}
	+\frac{2\ck_2}{\epsilon^2}\mynormb{\tilde{U}^{n,1}}_{L^2}
	\mynormb{\Delta\tilde{U}^{n,1}}_{L^2}\notag\\
	&\,+\frac{4\ck_5\ck_{\phi}\tau_n}{\epsilon^2}
	\sum_{j=2}^{k}\mynormb{\Delta\tilde{U}^{n,j}}_{L^2}
	\brab{\mynormb{\tilde{U}^{n,j}}_{L^2}+\mynormb{\tilde{U}^{n,j-1}}_{L^2}}	\notag\\
	&\,+\frac{\ck_6\tau_n}{\epsilon^2}\sum_{j=2}^{k}
	\mynormb{\Delta\tilde{U}^{n,j}}_{L^2}^2
	+\frac{\tau_n}{\lambda_{\min}\sigma_{\min}^{2}\epsilon^2}
	\sum_{j=1}^{k}\mynormb{\delta_{\tau}\zeta_{\mathrm{R}}^{n,j+1}}_{L^2}^2,
\end{align}
where the constant $\ck_6:=\frac{\ck_2^2}{\lambda_{\min}}$.
It is reasonable to assume further that the maximum time-step size is small such that $N\tau\le \ck_TT$. Also, let the constants $\ck_7:=
2^{2s-1}\exp\brab{2^{s-1}\ck_6T/\epsilon^2}\frac{\sqrt{T}s_{\mathrm{I}}\ck_1}{\sigma_{\min}\sqrt{\lambda_{\min}}}$, $\ck_8:=2^{2s+1}\exp\brab{2^{s-1}\ck_6T/\epsilon^2}\ck_4\ck_5\ck_{\phi}Ts_{\mathrm{I}}$ and $\ck_9:=4^{s}\exp\brab{2^{s-1}\ck_6T/\epsilon^2}\ck_2\ck_4\ck_TT.$
One can claim from the inequality \eqref{system: H2 time error R-IERK-difference2} that
\begin{align}\label{system: H2 time error R-IERK-difference3}
	\mynormb{\Delta\tilde{U}^{n,k+1}}_{L^2}\le&\,  (\ck_7\epsilon+\ck_8)\tau/\epsilon^2+\ck_9/\epsilon^2\quad\text{if $\tau\le \epsilon^2/\ck_6$.}
\end{align}
For any fixed $n$, two different cases are considered: (Case H2-1) if the $H_{per}^2$ semi-norm $$\mynormb{\Delta\tilde{U}^{n,j+1}}_{L^2}\le \frac{\sqrt{T}}{\sigma_{\min}\sqrt{\lambda_{\min}}\epsilon} \max\limits_{1\le j\le s_{\mathrm{I}}}\mynormb{\delta_{\tau}\zeta_{\mathrm{R}}^{n,j+1}}_{L^2}\quad\text{for $1\le j\le k$,}$$ the $H_{per}^2$ semi-norm  bound \eqref{system: H2 time error R-IERK-difference3} follows immediately. (Case H2-2) Otherwise, 
one can set $$\frac{\sigma_{\min}\sqrt{\lambda_{\min}}\epsilon}{\sqrt{T}} \mynormb{\Delta\tilde{U}^{n,j+1}}_{L^2}\ge \max\limits_{1\le j\le s_{\mathrm{I}}}\mynormb{\delta_{\tau}\zeta_{\mathrm{R}}^{n,j+1}}_{L^2}
\ge \mynormb{\delta_{\tau}\zeta_{\mathrm{R}}^{n,j+1}}_{L^2}\quad\text{for $1\le j\le k$,}$$
such that the inequality \eqref{system: H2 time error R-IERK-difference2} can be reduced into,
cf. the derivations from \eqref{system: L2 time error R-IERK-difference2} to \eqref{system: L2 time error R-IERK-difference3},
\begin{align*}
	\mynormb{\Delta\tilde{U}^{n,k+1}}_{L^2}\le&\,\mynormb{\Delta\tilde{U}^{n,1}}_{L^2}
	+\frac{\ck_6\tau_n}{\epsilon^2}\sum_{j=2}^{k}
	\mynormb{\Delta\tilde{U}^{n,j}}_{L^2}
	+\frac{\tau_n}{\sigma_{\min}\sqrt{\lambda_{\min}T}\epsilon}
	\sum_{j=1}^{k}\mynormb{\delta_{\tau}\zeta_{\mathrm{R}}^{n,j+1}}_{L^2}
	\notag\\
	&\,+\frac{4\ck_2}{\epsilon^2}\max_{1\le j\le k}\mynormb{\tilde{U}^{n,j}}_{L^2}	
	+\frac{4\ck_5\ck_{\phi}\tau_n}{\epsilon^2}\sum_{j=2}^{k}
	\brab{\mynormb{\tilde{U}^{n,j}}_{L^2}+\mynormb{\tilde{U}^{n,j-1}}_{L^2}}.
\end{align*}
By using the defect bound \eqref{system: rough assumption} 
and the $L^2$ norm estimate \eqref{system: L2 time error R-IERK-difference4}, the discrete Gr\"{o}nwall inequality in Lemma \ref{lemma: gronwall multi-stage RK1} 
together with the maximum time-step condition $\tau\le \epsilon^2/\ck_6$ gives
\begin{align*}
	\mynormb{\Delta\tilde{U}^{n,k+1}}_{L^2}\le&\, 4^{s-1}
	\exp\brab{2^{s-1}\ck_6T/\epsilon^2}
	\braB{\frac{2\sqrt{T}s_{\mathrm{I}}\ck_1\epsilon}{\sigma_{\min}\sqrt{\lambda_{\min}}}\tau
		+8\ck_4\ck_5\ck_{\phi}Ts_{\mathrm{I}}\tau
		+4\ck_2\ck_4n\tau}/\epsilon^2\\
		\leq&\,(\ck_7\epsilon+\ck_8)\tau/\epsilon^2+\ck_9/\epsilon^2.
\end{align*}
The above two cases (Case H2-1) and (Case H2-2) verify the $H_{per}^2$ semi-norm bound \eqref{system: H2 time error R-IERK-difference3}.

By combining  \eqref{system: L2 time error R-IERK-difference4} with \eqref{system: H2 time error R-IERK-difference3}, one can get the desired $H_{per}^2$ norm bound $\mynormb{\tilde{U}^{n,k+1}}_{H^2}\le \ck_{10}/\epsilon^2$ by taking $\ck_{10}:=\ck_4\epsilon^2/\ck_3+(\ck_7\epsilon+\ck_8)\epsilon^2/\ck_6+\ck_9$. That is, the estimate \eqref{bound: H2 time error hypothesis} holds for $l=n$ and $\ell=k+1$ and the mathematical induction arrives at the $H_{per}^2$ norm bound \eqref{bound: H2 time error}.

\subsubsection{$H_{per}^{m+4}$ norm boundedness of stage solutions}
It is to note that, the above proof of the $H_{per}^2$ norm bound \eqref{bound: H2 time error} can be extended to derive the following $H_{per}^3$ norm bound of stage solutions via the mathematical induction,
\begin{align*}
	\mynormb{\tilde{U}^{n,i}}_{H^3}\le \ck_{16}/\epsilon^2\quad\text{for $1\le n\le N$ and $1\le i\le s$.}
\end{align*}
Under the inductive hypothesis with $\ck_{16}^*:=\ck_{\phi}\epsilon^2+\ck_{\Omega}\ck_{16}$,
\begin{align*}
	\mynormb{\tilde{U}^{l,\ell}}_{H^3}\le \ck_{16}/\epsilon^2\quad\text{such that}\quad
	\mynormb{\nabla U^{l,\ell}}_{L^\infty}\le \ck_{16}^*/\epsilon^2
	\quad\text{for $1\le l\le n$ and $1\le \ell\le k$,}
\end{align*}
one can follow the proof of Lemma \ref{lemma: nonlinear term} to establish similar bounds for   
$\mynormb{\nabla F'(\Phi^{n,j})-\nabla F'(U^{n,j})}_{L^2}$ 
and $\mynormb{\delta_{\tau}\nabla F'(\Phi^{n,j})-\delta_{\tau}\nabla F'(U^{n,j})}_{L^2}$.
With these bounds of nonlinear term and the discrete Gr\"{o}nwall inequality in Lemma \ref{lemma: gronwall multi-stage RK1}, one can derive the $H^1$ semi-norm error bound $\mynormb{\nabla\tilde{U}^{n,k+1}}_{L^2}\le \ck_{12}\tau$ (under the maximum time-step size $\tau\le \epsilon^2/\ck_{11}$) by making the inner product of \eqref{system: time error R-IERK-difference} by $2\Delta\tilde{U}^{n,i+\frac{1}{2}}$ and following the derivations in \textbf{($L^2$ norm rough estimate)}. 
After that, we can make the inner product of \eqref{system: DOC time error R-IERK-difference} by $\frac{2}{\epsilon^2}\delta_{\tau}\Delta\tilde{U}^{n,j+1}$ and follow the derivations of  \textbf{($H^2$ seminorm rough estimate)} to get $H_{per}^3$ semi-norm error bound $\mynormb{\nabla\Delta\tilde{U}^{n,k+1}}_{L^2}\le \ck_{14}\tau/\epsilon^2+\ck_{15}/\epsilon^2$ (under the maximum time-step size $\tau\le \epsilon^2/\ck_{13}$ and $\tau\le \ck_TT/N$). Then we can complete the mathematical induction by taking the constant $\ck_{16}:=\ck_{12}\epsilon^2/\ck_{11}+\ck_{14}\epsilon^2/\ck_{13}+\ck_{15}$.

By performing the above process repeatedly, one can check that there exists a positive constant $\Ck_{\phi}$, independent of  $\tau_n$, such that the stage solutions $U^{n,i}$ fulfill $\mynormb{U^{n,i}}_{H^{m+4}}\le \Ck_{\phi}/\epsilon^2$ for $1\le n\le N$ and $2\le i\le s$. Combining it with the time-discrete system \eqref{system: time discrete R-IERK}, we have the following theorem. 
\begin{theorem}\label{thm: regularity of time-discrete solution}
	Assume that the solution of the CH problem \eqref{cont: Problem-CH} satisfies the regularity assumption \eqref{regular: exact solution} with $m\ge1$, and the differentiation matrix $D_{\mathrm{R}}=A_{\mathrm{E}}^{-1}E_{s_{\mathrm{I}}}$ of the R-IERK method \eqref{Scheme: R-IERK stabilized} and the inverse $D_{\mathrm{R}}^{-1}$ are positive definite. If the maximum time-step size $\tau$ is sufficiently small, there exists a positive constant $\Ck_{\phi}$, independent of the time-step sizes $\tau_n$, such that the stage solutions $U^{n,i}$ of the time-discrete system \eqref{system: time discrete R-IERK} are bounded,
	\begin{align*}
		\mynormb{U^{n,i}}_{H^{m+4}}+\mynormb{(U^{n,i}-U^{n,1})/\tau_n}_{H^{m}}\le \Ck_{\phi}/\epsilon^2
		\quad\text{for $1\le n\le N$ and $2\le i\le s$.}
	\end{align*}	
	\end{theorem}

\subsection{Uniform boundedness of stage solutions and discrete energy law}

\subsubsection{Fully discrete error system}
As seen from Lemma \ref{lemma: IERK energy stability},  proving the uniform boundedness of stage solutions $u_h^{n,i}$ is necessary to establish the energy dissipation laws for the fully-discrete R-IERK method \eqref{Scheme: R-IERK stabilized}, which can be viewed as the spatial approximation of the time discrete system \eqref{system: time approximation R-IERK}.

In general, we evaluate the stage error  by a usual splitting,
$U^{n,\ell}-u_h^{n,\ell}=U^{n,\ell}-U_M^{n,\ell}+e_h^{n,\ell},$
where $U_M^{n,\ell}:=P_MU^{n,\ell}$ is the $L^2$ projection of time-discrete solution $U^{n,\ell}$ and  $e_h^{n,\ell}:=U_M^{n,\ell}-u_h^{n,\ell}\in \mathbb{\mathring V}_{h}$
is the difference between the projection $U_M^{n,\ell}$
and the solution $u_h^{n,\ell}$ of the R-IERK method 
\eqref{Scheme: R-IERK stabilized}. Without losing the generality, we set the initial data $u_h^{1,1}:=U_M^{1,1}$ such that $e_h^{1,1}=0.$
Note that,  the $L^2$ projection solution $U_M^{n,\ell}\in\mathscr{F}_M$
and  the volume conservation \eqref{eq: stage volume preserving}	arrive at
\[
\myinnerb{U_M^{n,\ell},1}
=\myinnerb{U_M^{1,1},1}=\myinnerb{u^{1,1},1}
=\myinnerb{u^{n,\ell},1},
\]
so that the error function $e^{n,\ell}\in \mathbb{\mathring V}_{h}$.
Applying Lemma \ref{lem:Projection-Estimate} with 
$U^{n,\ell}\in C^{2}\brab{[0,T];{H}_{per}^{m}}$, one has
$$\mynormb{U^{n,\ell}-U_M^{n,\ell}}=\mynormb{I_M\brat{U^{n,\ell}-U_M^{n,\ell}}}_{L^2}
\le \Ck_{\phi}\mynormb{I_MU^{n,\ell}-U_M^{n,\ell}}_{L^2}\le \Ck_{\phi}h^{m}\absb{U^{n,\ell}}_{H^{m}}.$$
Once the upper bound of $\mynorm{e^n}$ is available, the $L^2$ norm error estimate follows immediately,
\begin{align}\label{Triangle-Projection-Estimate}
	\mynormb{U^{n,\ell}-u^{n,\ell}}\le&\,
	\ck_{\phi}h^{m}/\epsilon^2+\mynormb{e^{n,\ell}}\quad\text{for $1\le n\le N$ and $1\le \ell\le s_{\mathrm{I}}$.}
\end{align}

 To bound the $L^2$ norm of $e_h^{n,\ell}$, we consider the space consistency error for a semi-discrete system having a projected solution $U_M$.
A substitution of the $L^2$ projection solution $U_M$ and differentiation  operator $\Delta_h$ into the equation \eqref{system: time discrete R-IERK} yields the discrete system
\begin{align}\label{system: projection time discrete R-IERK}				
	U_M^{n,i+1}=&\,U_M^{n,1}+\tau_n\sum_{j=1}^{i}\hat{a}_{i+1,j}\Delta_h \kbraB{L_{\kappa,h}U_M^{n,j+\frac1{2}}-f_{\kappa}(U_M^{n,j})}+\tau_n\zeta_P^{n,i+1}
\end{align}
for $n\ge1$ and $1\le i\le s_{\mathrm{I}}$, where $\zeta_{P}^{n,i+1}=\zeta_{P}(\myvec{x}_h,t_{n,i+1})$ represents the spatial consistency error arising from the $L^2$ projection of time-discrete solution $U^{n,j}$, that is,
\begin{align*}
	\zeta_{P}^{n,i+1}:=&\,\frac1{\tau_n}(U_M^{n,i+1}-U_M^{n,1})
	-\sum_{j=1}^{i}\hat{a}_{i+1,j}\Delta_h \kbraB{L_{\kappa,h}U_M^{n,j+\frac1{2}}-f_{\kappa}(U_M^{n,j})}\\
	&\,-\frac1{\tau_n}(U^{n,i+1}-U^{n,1})+\sum_{j=1}^{i}\hat{a}_{i+1,j}\Delta \kbraB{L_{\kappa}U^{n,j+\frac1{2}}-f_{\kappa}(U^{n,j})}\quad\text{for $\myvec{x}_h\in\Omega_{h}$.}
\end{align*}
By Theorem \ref{thm: regularity of time-discrete solution}
and Lemma \ref{lem:Projection-Estimate}, one can follow the proof of \cite[Theorem 3.1]{LiaoJiZhang:2022PFC} to obtain
\begin{align}\label{Projection-truncation error}
	\mynormb{\zeta_{P}^{n,i+1}}\le \Ck_\phi h^m\max_{1\le i\le s_{\mathrm{I}}}\brab{\mynormb{U^{n,i+1}}_{H^{m+4}}+\tfrac1{\tau_n}\mynormb{U^{n,i+1}-U^{n,1}}_{H^{m}}}\le \hck_{1} h^m/\epsilon^2.
\end{align}

 Subtracting the fully discrete scheme \eqref{Scheme: R-IERK stabilized} from \eqref{system: projection time discrete R-IERK}, one has the following error system 
\begin{align*}				
	e_h^{n,i+1}=&\,e_h^{n,1}+\tau_n\sum_{j=1}^{i}\hat{a}_{i+1,j}\Delta_h \kbraB{L_{\kappa,h}e_h^{n,j+\frac1{2}}+f_{\kappa}(u_h^{n,j})-f_{\kappa}(U_M^{n,j})}+\tau_n\zeta_{P}^{n,i+1}
\end{align*} 
or, with the matrix $(\underline{\hat{a}}_{i+1,j})_{s_{\mathrm{I}}\times s_{\mathrm{I}}}=E_{s_{\mathrm{I}}}^{-1}A_{\mathrm{E}}=D_{\mathrm{R}}^{-1}$,
\begin{align}\label{Scheme: Error R-IERK stabilized}				
	\delta_{\tau}e_h^{n,i+1}=\tau_n\sum_{j=1}^{i}\underline{\hat{a}}_{i+1,j}\Delta_h \kbraB{L_{\kappa,h}e_h^{n,j+\frac1{2}}
		+f_{\kappa}(u_h^{n,j})-f_{\kappa}(U_M^{n,j})}
	+\tau_n\delta_{\tau}{\zeta_{P}}^{n,i+1}
\end{align} 
for $n\ge1$ and $1\le i\le s_{\mathrm{I}}$.

\subsubsection{Uniform boundedness of fully discrete solutions}

We use the complete mathematical induction to prove that the stage solutions $u_h^{n,i}$ of the fully discrete scheme \eqref{Scheme: R-IERK stabilized} are bounded in the $L^2$ norm, that is,
\begin{align}\label{bound: L2 error rough}
	\mynormb{e^{n,i}}\le \hck_{4}h^m/\epsilon^2\quad\text{for $1\le n\le N$ and $1\le i\le s$.}
\end{align}
Obviously, it holds for $n=1$ and $i=1$ since $e^{1,1}=0$. Put the inductive hypothesis
\begin{align}\label{bound: L2 error rough hypothesis}
	\mynormb{e^{l,\ell}}\le \hck_{4}h^m/\epsilon^2
	\quad\text{for $1\le l\le n$ and $1\le \ell\le k$,}
\end{align}
such that $\mynormb{e^{l,\ell}}_{\infty}\le h^{-1}\mynormb{e^{l,\ell}}\le \hck_{4}h^{m-1}/\epsilon^2\le1$ if $h\le \sqrt[m-1]{\epsilon^2/\hck_4}$ and thus
\begin{align}\label{bound: maximum error rough hypothesis}
	\mynormb{u^{l,\ell}}_{\infty}\le \mynormb{U_M^{l,\ell}}_{\infty}+\mynormb{e^{l,\ell}}_{\infty}\le \hck_{4}^*/\epsilon^2
	\quad\text{for $1\le l\le n$ and $1\le \ell\le k$,}
\end{align}
where the constant $\hck_{4}^*:=\ck_{\phi}+\epsilon^2.$ In the following, we are to prove that $\mynormb{e^{n,k+1}}\le \hck_{4}h^m/\epsilon^2$.

Making the inner product of \eqref{Scheme: Error R-IERK stabilized} by $2e_h^{n,i+\frac{1}{2}}$, applying the discrete Green's formula and summing the stage index $i$ from $i=1$ to $k$, one can find that
\begin{align}\label{Scheme: L2 Error R-IERK stabilized}
	\mynormb{e^{n,k+1}}^2-\mynormb{e^{n,1}}^2
	=&\,2\tau_n\sum_{i=1}^{k}\sum_{j=1}^{i}\underline{\hat{a}}_{i+1,j}
	\myinnerb{
		L_{\kappa,h}e^{n,j+\frac{1}{2}},\Delta_h e^{n,i+\frac{1}{2}}}
	+2\tau_n\sum_{i=1}^{k}
	\myinnerb{\delta_{\tau}\zeta_{P}^{n,i+1},e^{n,i+\frac{1}{2}}}\notag\\
	&\,+2\tau_n\sum_{i=1}^{k}\sum_{j=1}^{i}
	\underline{\hat{a}}_{i+1,j}\myinnerb{f_{\kappa}(u^{n,j})-f_{\kappa}(U_M^{n,j}), \Delta_he^{n,i+\frac{1}{2}}}.
\end{align}
Lemma \ref{lemma: convolution terms} (i) yields
\begin{align*}
	2\tau_n\sum_{i=1}^{k}\sum_{j=1}^{i}\underline{\hat{a}}_{i+1,j}\myinnerb{
		L_{\kappa,h}e^{n,j+\frac{1}{2}},\Delta_he^{n,i+\frac{1}{2}}}
	\le-2\epsilon^2\sigma_{\min}\tau_n \sum_{i=1}^{k}\mynormb{\Delta_he^{n,i+\frac{1}{2}}}^2.
\end{align*}
For the last term in \eqref{Scheme: L2 Error R-IERK stabilized}, one can apply Lemma \ref{lemma: convolution terms} (ii) and Lemma \ref{lemma: nonlinear term} (i) to derive that
\begin{align*}
	2\tau_n\sum_{i=1}^{k}\sum_{j=1}^{i}\underline{\hat{a}}_{i+1,j}\myinnerb{
		f_{\kappa}(u^{n,j})-f_{\kappa}(U_M^{n,j}),\Delta_he^{n,i+\frac{1}{2}}}
	\le \frac{2\tau_n}{\lambda_{\min}} \sum_{i=1}^{k}\mynormb{f_{\kappa}(u^{n,i})-f_{\kappa}(U_M^{n,i})}
	\mynormb{\Delta_he^{n,i+\frac{1}{2}}}\\
	\le 
	2\epsilon^2\sigma_{\min}\tau_n
	\sum_{i=1}^{k}\mynormb{\Delta_he^{n,i+\frac{1}{2}}}^2
	+\frac{\tau_n}{2\sigma_{\min}\lambda_{\min}^2\epsilon^2}	 \sum_{i=1}^{k}\mynormb{e^{n,i}}^2
	\int_0^1\brab{\kappa+\mynormb{F''\brab{\mathrm{Q}_{\gamma,h}^{n,i}[u]}}_{\infty}}^2\zd\gamma,
\end{align*}
where the discrete functional $\mathrm{Q}_{\gamma,h}^{n,i}[u]$ is defined via Lemma \ref{lemma: nonlinear term}, that is, $\mathrm{Q}_{\gamma,h}^{n,j}[u]:=U_M^{n,j}+\gamma e_h^{n,j}$, so that $\mynormb{F''\brab{\mathrm{Q}_{\gamma,h}^{n,i}[u]}}_{\infty}\le \hck_2:=\max\limits_{\mynormt{\xi}_{\infty}\le \hck_{4}^*/\epsilon^2}\mynormb{F''(\xi)}_{\infty}$ according to the maximum norm bound \eqref{bound: maximum error rough hypothesis}.
Thus, it follows from \eqref{Scheme: L2 Error R-IERK stabilized} that
\begin{align}\label{Scheme: L2 Error R-IERK stabilized2}
	\mynormb{e^{n,k+1}}^2-\mynormb{e^{n,1}}^2
	\le\frac{\hck_3\tau_n}{\epsilon^2} \sum_{i=1}^{k}\mynormb{e^{n,i}}^2
	+2\tau_n\sum_{i=1}^{\ell}\mynormb{\delta_{\tau}\zeta_{P}^{n,i+1}}\mynormb{e^{n,i+1}},
\end{align}
where the constant $\hck_3:=\frac{(\hck_2+\kappa)^2}{2\sigma_{\min}\lambda_{\min}^2}$. For any index $n$, choose $k_0$ satisfying $1\le k_0\le k$ such that $\mynormb{e^{n,k_0+1}}:=\max_{0\le \ell\le k}\mynormb{e^{n,\ell+1}}$. We set $k=k_0$ in \eqref{Scheme: L2 Error R-IERK stabilized2} to get
\begin{align*}
	\mynormb{e^{n,k_0+1}}\le\mynormb{e^{n,1}}
	+\frac{\hck_3\tau_n}{\epsilon^2} \sum_{i=1}^{k_0}\mynormb{e^{n,i}}
	+2\tau_n\sum_{i=1}^{k_0}\mynormb{\delta_{\tau}\zeta_{P}^{n,i+1}},
\end{align*}
and then, due to $k_0\le k$ and $\mynormb{e^{n,k+1}}\le\mynormb{e^{n,k_0+1}}$,
\begin{align*}
	\mynormb{e^{n,k+1}}\le\mynormb{e^{n,1}}
	+\frac{\hck_3\tau_n}{\epsilon^2}\sum_{i=1}^{k}\mynormb{e^{n,i}}
	+2\tau_n\sum_{i=1}^{k}\mynormb{\delta_{\tau}\zeta_{\mathrm{P}}^{n,i+1}}.
\end{align*}
By using the truncation error bound \eqref{Projection-truncation error}, the discrete Gr\"{o}nwall inequality in Lemma \ref{lemma: gronwall multi-stage RK1} 
together with the maximum time-step condition $\tau\le \epsilon^2/\hck_3$ gives
\begin{align*}
	\mynormb{e^{n,k+1}}\le&\, 4^{s}
	\exp\brab{2^{s-1}\hck_3t_{n-1}/\epsilon^2}s_{\mathrm{I}}\hck_1t_nh^m/\epsilon^2\le \hck_4h^m/\epsilon^2,
\end{align*}
where the constant $\hck_4:=4^{s}
\exp\brab{2^{s-1}\hck_3T/\epsilon^2}s_{\mathrm{I}}\hck_1T$.

It says that the estimate \eqref{bound: L2 error rough hypothesis} holds for $l=n$ and $\ell=k+1$. Thus the principle of mathematical induction confirms the $L^2$ norm error estimate \eqref{bound: L2 error rough}. The inverse estimate yields 
\begin{align}\label{bound: maximum error rough}
	\mynormb{u^{n,i}}_{\infty}\le \hck_{4}^*/\epsilon^2
	\quad\text{for $1\le n\le N$ and $1\le i\le s$}\quad\text{if $h\le \sqrt[m-1]{\epsilon^2/\hck_4}$.}
\end{align}


Thanks to the maximum norm bound \eqref{bound: maximum error rough} of stage solutions,
one can take the positive constant $\ck_0:=\hck_{4}^*/\epsilon^2$ in Lemma \ref{lemma: IERK energy stability} 
to obtain the following result.
\begin{theorem}\label{thm: RIERK energy stability}
	Assume that the solution of the CH problem \eqref{cont: Problem-CH} satisfies the regularity assumption \eqref{regular: exact solution} with $m\ge1$, and the differentiation matrix $D_{\mathrm{R}}=A_{\mathrm{E}}^{-1}E_{s_{\mathrm{I}}}$ of the R-IERK method \eqref{Scheme: R-IERK stabilized} and the inverse $D_{\mathrm{R}}^{-1}$ are positive definite.
	If the spatial length $h$  and the maximum step size $\tau$ are sufficiently small, and the parameter $\kappa$ in \eqref{def: stabilized parameter} is chosen properly large such that 
	$\kappa\ge\max_{\mynormt{\xi}_{\infty}\le \hck_{4}^*/\epsilon^2}\mynormb{F''(\xi)}_{\infty},$
	then the stage solutions $u_h^{n,i}$ of the R-IERK method \eqref{Scheme: R-IERK stabilized} are bounded in the maximum norm and they preserve 
	the original energy dissipation law \eqref{cont:energy dissipation}
	at all stages,
	\begin{align*}		
		E[u^{n,j+1}]-E[u^{n,1}]\le&\,
		\frac1{\tau}\sum_{k=1}^{j}\myinnerB{\Delta_h^{-1}\delta_{\tau}u^{n,k+1},
			\sum_{\ell=1}^{k}d_{k\ell}^{(R)}\delta_{\tau}u^{n,\ell+1}}
		\quad\text{for $1\le n\le N$ and $1\le j\le s_{\mathrm{I}}$.}
	\end{align*}
\end{theorem}

\subsection{$L^2$ norm estimate with full accuracy}

For the convergence of the R-IERK method \eqref{Scheme: R-IERK stabilized}, we will overlook the defects at internal stages as done in the literature, cf. \cite{DuJuLiQiao:2021SIREV,FuTangYang:2024}, and assume that the stage temporal defects $\zeta_{\mathrm{R}}^{n,i+1}$ satisfy 
\begin{align}\label{Estimate: refined assumption}
	\zeta_{\mathrm{R}}^{n,i+1}=0\quad\text{for $1\le i\le s_{\mathrm{I}}$}\quad \text{and}\quad \mynormb{\zeta_{\mathrm{R}}^{n,s}}\le \ck_2\tau ^p
	\quad \text{for $1\le n\le N$.}
\end{align}
Actually, the rough setting \eqref{system: rough assumption} is not enough to derive the sharp error estimate at time levels with fully accuracy.
With the help of Theorem \ref{thm: RIERK energy stability}, the following result can be verified by following the $L^2$ norm estimate in subsection 4.3. The technical details are left to interested readers.

\begin{theorem}\label{thm: RIERK convergence}
	Assume that the solution of the CH problem \eqref{cont: Problem-CH} satisfies the regularity assumption \eqref{regular: exact solution} with $m\ge1$, and the differentiation matrix $D_{\mathrm{R}}=A_{\mathrm{E}}^{-1}E_{s_{\mathrm{I}}}$ of the R-IERK method \eqref{Scheme: R-IERK stabilized} and the inverse $D_{\mathrm{R}}^{-1}$ are positive definite.
	If the spatial length $h$  and the maximum step size $\tau$ are sufficiently small, and the parameter $\kappa$ in \eqref{def: stabilized parameter} is chosen properly such that 
	$\kappa\ge\max_{\mynormt{\xi}_{\infty}\le \hck_{4}^*/\epsilon^2}\mynormb{F''(\xi)}_{\infty},$
	then the solution $u_h^n$ of R-IERK method \eqref{Scheme: R-IERK stabilized} is convergent in the $L^2$ norm with an order of $\mathcal{O}(\tau^p+h^m)$.
\end{theorem}

\section{Numerical experiments}
\label{sec: experiments}
\setcounter{equation}{0}


\begin{example}\label{1ex: FuShenYang}
	Consider the Cahn-Hilliard model \eqref{cont: Problem-CH} with an exterior force $g(x,y,t)$ subject to the initial data $\Phi^{0} = \sin (\pi x) \sin (\pi y)$ on $\Omega=(0,2)^2$ with the interface parameter $\epsilon=0.2$. The source term $g$ is set by choosing the exact solution $\Phi(x,y;t) = e^{-t}\sin (\pi x) \sin (\pi y)$. Always, the spatial operators are approximated by the Fourier pseudo-spectral approximation with $64\times64$ grid points. 
\end{example}

\begin{figure}[htb!]
	\centering
	\subfigure[R-IERK(2,4;$c_2$) method  \eqref{scheme: R-IERK-2-4-c2}]{
		\includegraphics[width=2.3in]{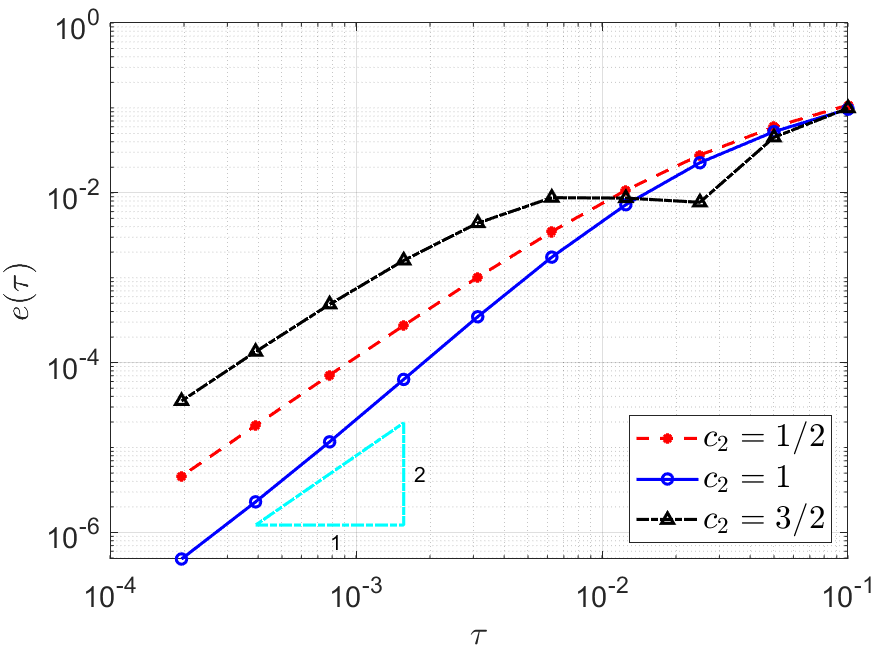}}
	\subfigure[R-IERK(3,6;$\hat{a}_{52}$) method  \eqref{scheme: R-IERK-3-6-ah52}]{
		\includegraphics[width=2.3in]{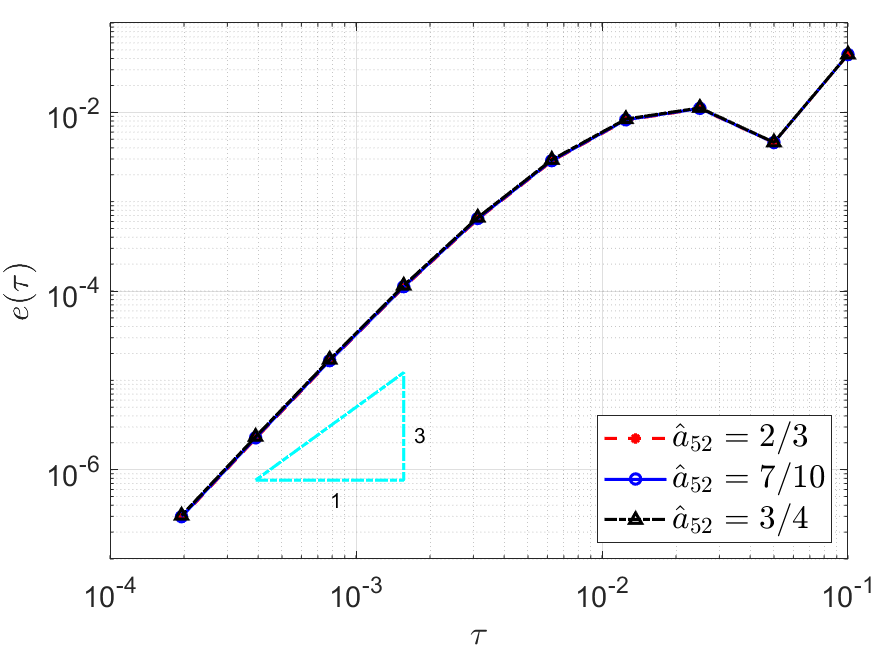}}
	\caption{Solution errors of the two R-IERK methods with different parameters.}
	\label{fig: IERK2 convergence}
\end{figure}

We examine the convergence of the R-IERK(2,4;$c_2$) method \eqref{scheme: R-IERK-2-4-c2} and the R-IERK(3,6;$\hat{a}_{52}$) method \eqref{scheme: R-IERK-3-6-ah52} by choosing the final time $T=1$ and the stabilized parameter $\kappa=4$. Figure \ref{fig: IERK2 convergence} lists
the $L^\infty$ norm error $e(\tau):=\max_{1\le n\le N}\mynorm{\Phi_h^{n} - \Phi(t_n)}_\infty$ for the two R-IERK methods on
halving time steps $\tau=2^{-k}/10$ for $0\le k\le 9$. As expected, the R-IERK(2,4;$c_2$) method \eqref{scheme: R-IERK-2-4-c2} and R-IERK(3,6;$\hat{a}_{52}$) method \eqref{scheme: R-IERK-3-6-ah52} are second-order and third-order accurate in time, respectively. It suggests that the different
parameters for the R-IERK(2,4;$c_2$) method \eqref{scheme: R-IERK-2-4-c2} would arrive at different precision and the case $c_2=1$ generates the most accurate solution when $\tau$ is smaller than $10^{-2}$; while the R-IERK(3,6;$\hat{a}_{52}$) method \eqref{scheme: R-IERK-3-6-ah52} with different parameters generates almost the same solution.

\subsection{Tests of R-IERK(2,4;$c_2$) methods}

\begin{figure}[htb!]
	\centering
	\subfigure[discrete energy]{
		\includegraphics[width=2.3in]{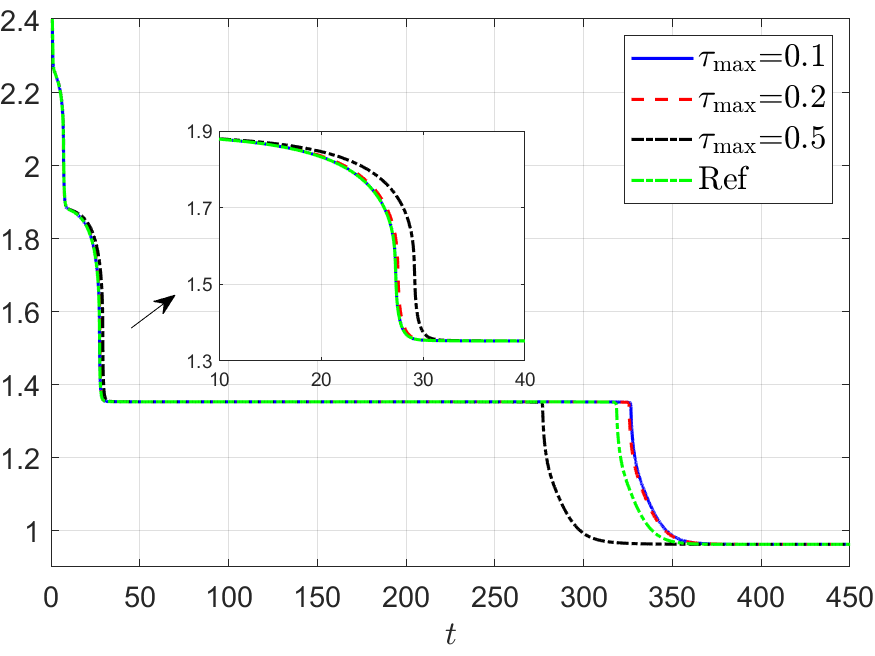}}\hspace{10mm}
	\subfigure[adaptive step sizes]{
		\includegraphics[width=2.3in]{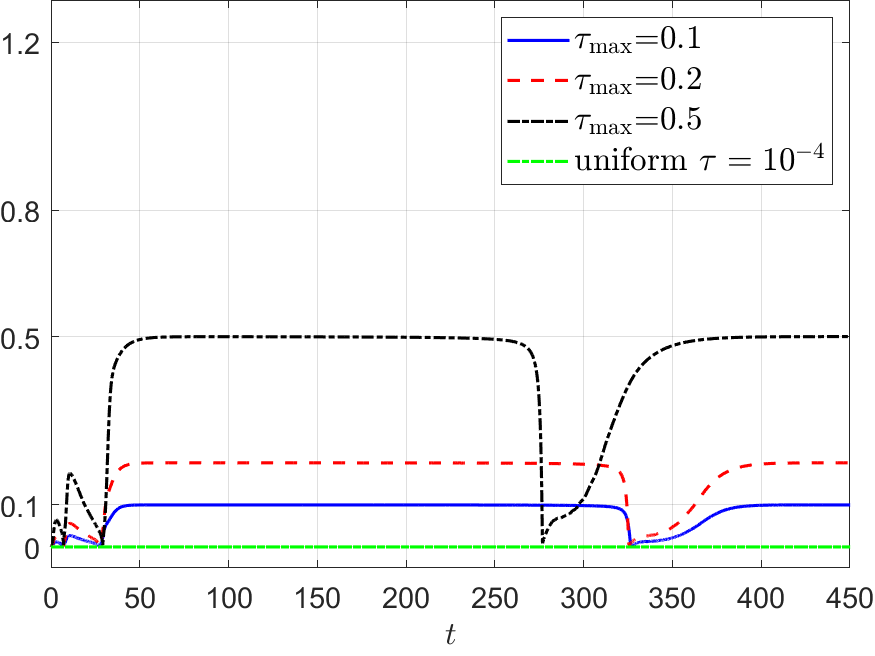}}
	\caption{ Energy curves and adaptive time-step sizes of R-IERK(2,4;1) method.}
	\label{fig: IERK2-I decay a33}
\end{figure}

We examine the discrete energy behaviors by running the IERK(2,3) method and the R-IERK(2,4;$c_2$) method \eqref{scheme: R-IERK-2-4-c2} with three different parameters $c_2=\frac12,1$ and $\frac32$ by the same adaptive time-stepping strategy adopted in Example \ref{1example: energy_2}.
The reference solution is generated by the IERK(2,3) method with $\tau=10^{-4}$, while in the adaptive time-stepping algorithm, we consider the three different scenes: (a) $\tau_{\max}=0.1$, (b) $\tau_{\max}=0.2$ and (c) $\tau_{\max}=0.5$.

At first, we run a special R-IERK(2,4;$c_2$) method (taking $c_2=1$ for example) for the three different scenes. As expected, see Figure \ref{fig: IERK2-I decay a33}, the discrete energy curves are decreasing, and the energy curve is closer to the reference energy when a smaller $\tau_{\max}$ is adopted. The CPU time and the total number of time-level listed in Table \ref{tabel: IERK2} show the effectiveness and efficiency of the R-IERK(2,4;1) method \eqref{scheme: R-IERK-2-4-c2} with $T=450$. Surprisingly, the R-IERK(2,4;1) method preforms well and generates reliable energy curves although it can not be covered by our theory, see Remark \ref{remark: zero-eigenvalue c2_1}.

\begin{table}[htb!]
	\centering
	\caption{CPU time and total time-level of R-IERK(2,4;1) method.}
	\begin{tabular}{@{}lcccc@{}}
		\toprule
		Time-stepping strategy        & $\tau = 10^{-4}$ & $\tau_{\max} = 0.1$ & $\tau_{\max} = 0.2$ & $\tau_{\max} = 0.5$ \\ \midrule
		CPU time (in seconds)          & 6724.63           & 11.90     & 6.42       &   2.62     \\[2pt]
		Time levels        & $4.5\times10^6$             & 11714        & 5840         & 2304           \\ \bottomrule
	\end{tabular}
	\label{tabel: IERK2}
\end{table}

\begin{figure}[htb!]
	\centering
	\subfigure[$\tau_{\max}=0.1$]{
		\includegraphics[width=2.05in]{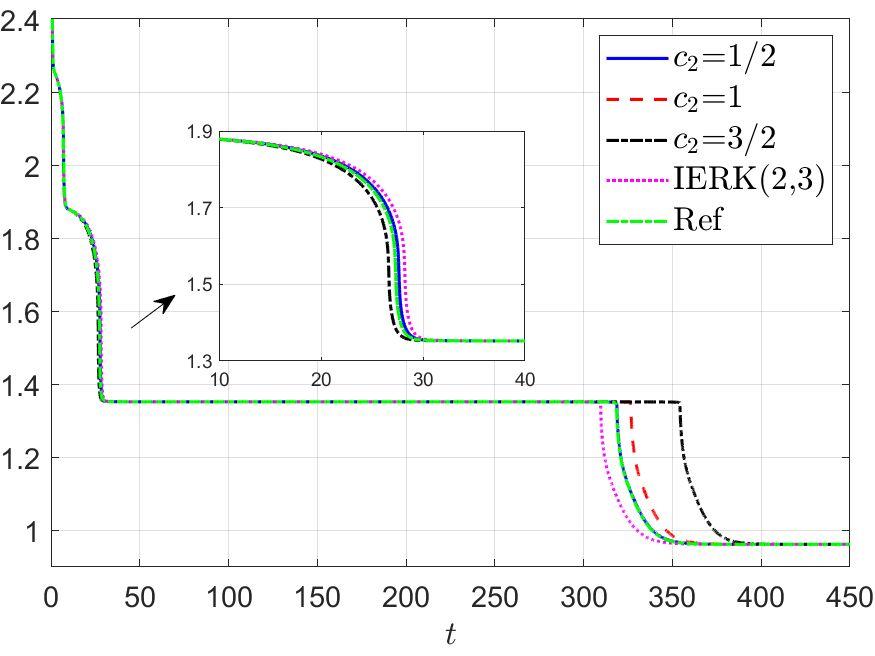}}
	\subfigure[$\tau_{\max}=0.2$]{
		\includegraphics[width=2.05in]{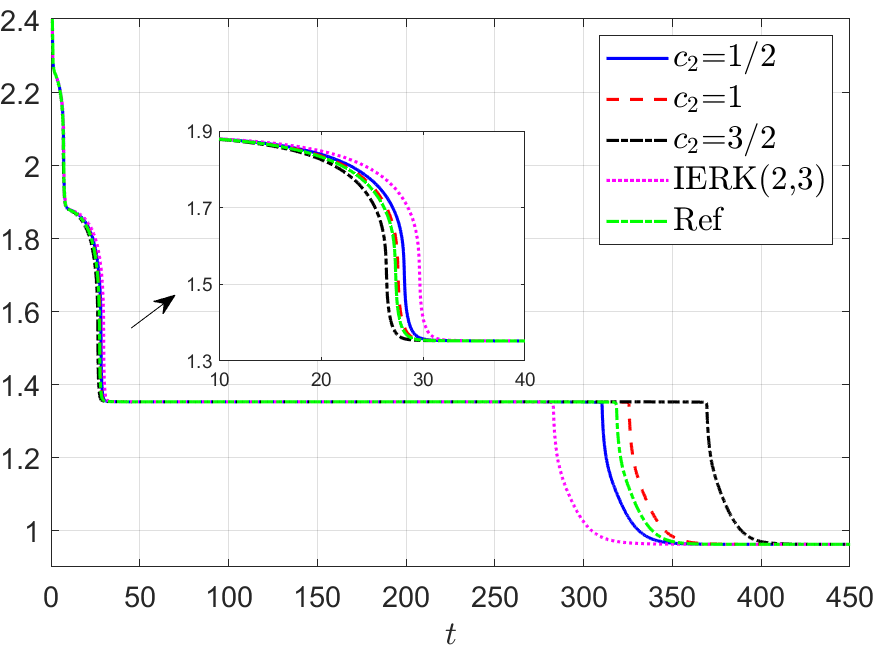}}
	\subfigure[$\tau_{\max}=0.5$]{
		\includegraphics[width=2.05in]{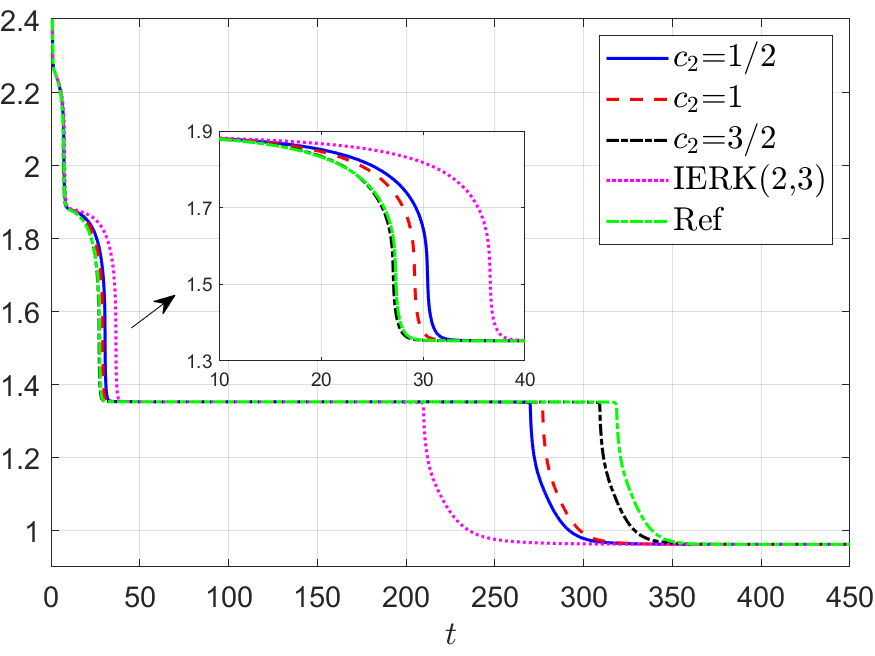}}
	\caption{Energy behaviors of R-IERK(2,4;$c_2$) method \eqref{scheme: R-IERK-2-4-c2} and IERK(2,3) method.}
	\label{1fig: IERK decay c2}
\end{figure}

Figure \ref{1fig: IERK decay c2} depicts the energy curves generated by the IERK(2,3) method and the R-IERK(2,4;$c_2$) method \eqref{scheme: R-IERK-2-4-c2} with different parameters $c_2=\tfrac12, 1$ and $\tfrac32$ for three different maximum time-steps: (a) $\tau_{\max}=0.1$, (b) $\tau_{\max}=0.2$ and (c) $\tau_{\max}=0.5$. As expected, 
the discrete energy curves are all decreasing. 
It is seen that the energy curves of the IERK(2,3) method have significant changes for different maximum time-steps, while the discrete energy curves computed by the R-IERK(2,4;$c_2$) method \eqref{scheme: R-IERK-2-4-c2} are relatively robust with respect to the change of time-steps. They suggest that the R-IERK(2,4;$c_2$) methods \eqref{scheme: R-IERK-2-4-c2} with the  $\tau_n\lambda_{\mathrm{ML}}$-independent average dissipation rates allow some larger adaptive step-sizes and would be more preferable in adaptive time-stepping simulations.



\subsection{Tests of R-IERK(3,6;$\hat{a}_{52}$) methods}

We examine the energy behaviors of the R-IERK(3,6;$\hat{a}_{52}$) method \eqref{scheme: R-IERK-3-6-ah52} with three different parameters $\hat{a}_{52}=\tfrac23$, $\tfrac7{10}$ and $\tfrac34$ for Example \ref{1example: energy_2} with the adaptivity parameter $\eta=500$.
The reference solution is generated with $\tau=10^{-4}$ by the Lobatto-type IERK3-2 method in \cite{LiaoWangWen:2024arxiv} with the parameter $a_{43}=-\tfrac{3}5$, called IERK(3,5) method hereafter, which is regarded as the best one among the third-order IERK methods in \cite{LiaoWangWen:2024arxiv} with the associated average dissipation rate $\mathcal{R}^{(3,5)}_{\mathrm{L}}=\tfrac{5}{4} + \tfrac{2}{5}\tau\overline{\lambda}_{\mathrm{ML}}$. 

\begin{figure}[htb!]
	\centering
	\subfigure[$\tau_{\max}=0.2$]{
		\includegraphics[width=2.05in]{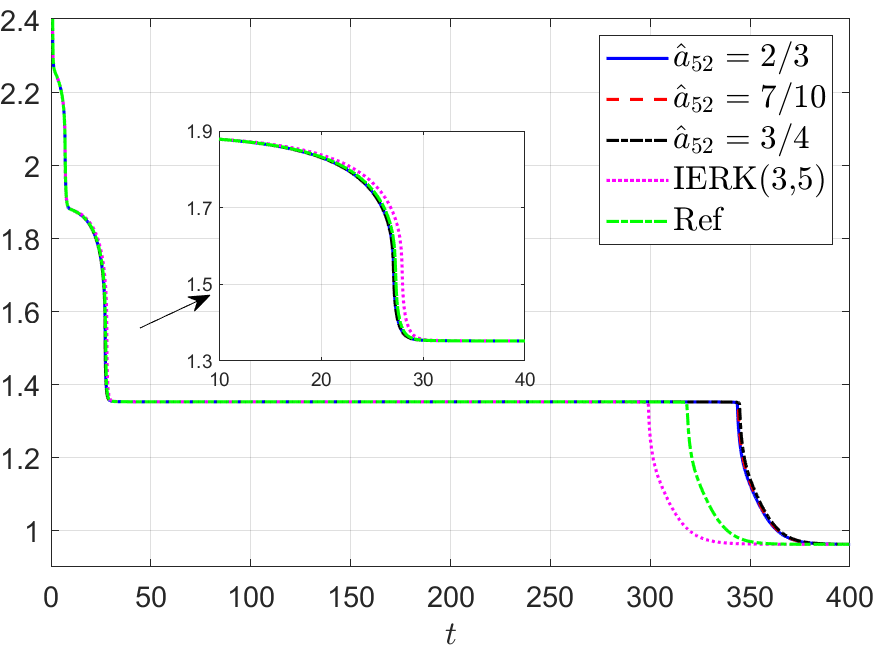}}
	\subfigure[$\tau_{\max}=0.5$]{
		\includegraphics[width=2.05in]{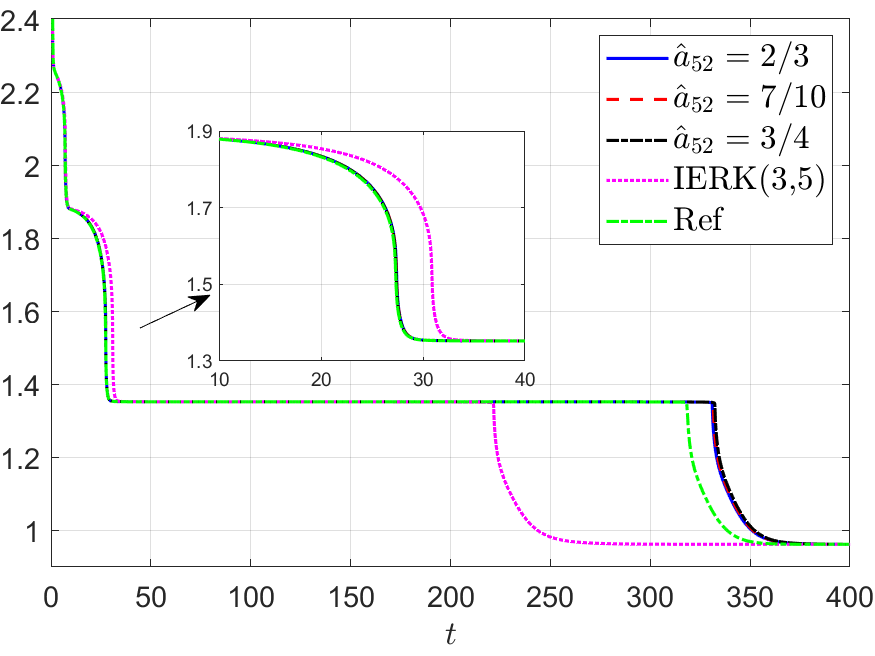}}
	\subfigure[$\tau_{\max}=0.8$]{
		\includegraphics[width=2.05in]{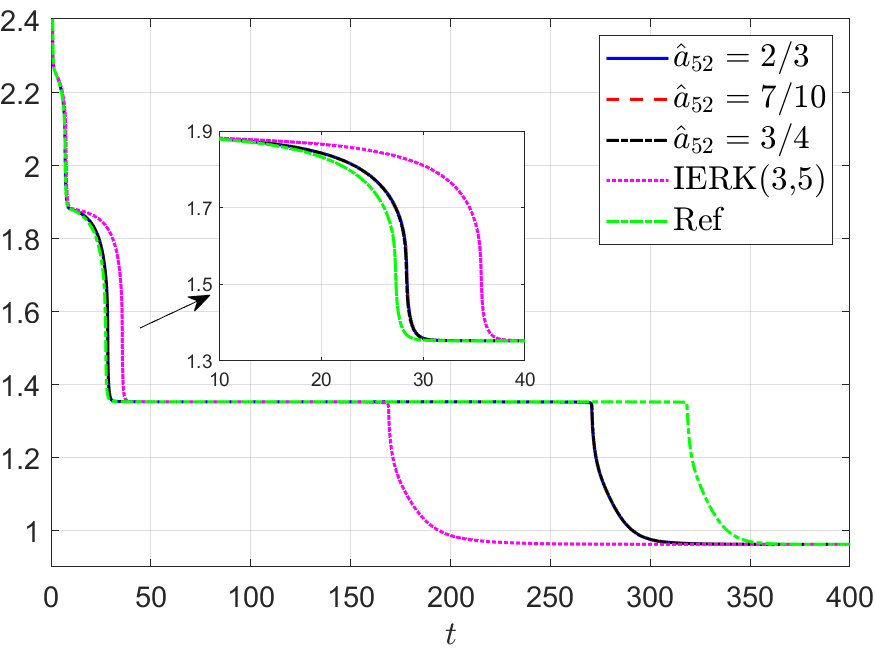}}
	\caption{Energy behaviors of R-IERK(3,6;$\hat{a}_{52}$) method \eqref{scheme: R-IERK-3-6-ah52} and IERK(3,5) method.}
	\label{fig: IERK decay a52}
\end{figure}

Figure \ref{fig: IERK decay a52} depicts the energy curves of the IERK(3,5) method and the R-IERK(3,6;$\hat{a}_{52}$) method \eqref{scheme: R-IERK-3-6-ah52} for three different scenes: (a) $\tau_{\max}=0.2$, (b) $\tau_{\max}=0.5$ and (c) $\tau_{\max}=0.8$.
As seen, the energy curves generated by the R-IERK(3,6;$\hat{a}_{52}$) method are indistinguishable for all scenes, which seems to be accordant with the convergence tests in Figure \ref{fig: IERK2 convergence} (b).
More importantly, one can observe that the energy curves of the IERK(3,5) method have significant changes for different maximum time-steps, while the discrete energy curves computed by the R-IERK(3,6;$\hat{a}_{52}$) methods are relatively robust with respect to the change of time-steps. They suggest that the R-IERK(3,6;$\hat{a}_{52}$) methods allow some larger adaptive step-sizes and would be more preferable in adaptive simulations.

\begin{figure}[htb!]
	\centering
	\subfigure{
		\includegraphics[width=1.8in]{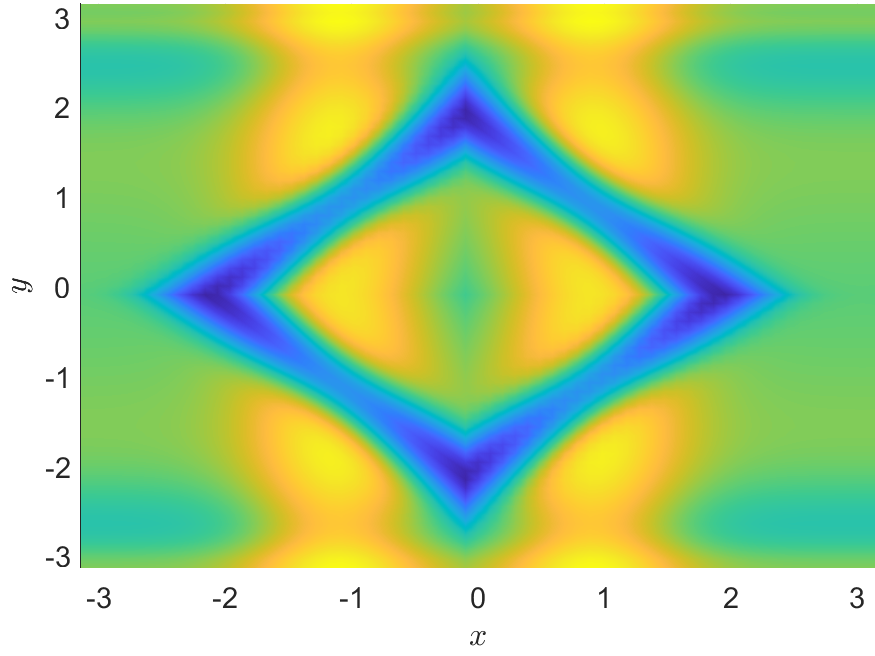}}
	\subfigure{
		\includegraphics[width=1.8in]{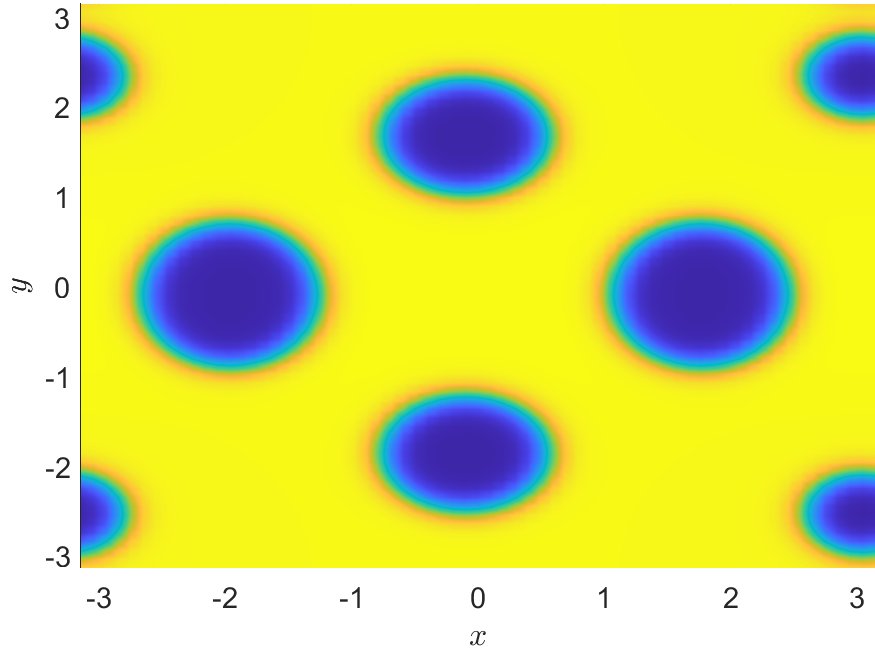}}
	\subfigure{
		\includegraphics[width=1.8in]{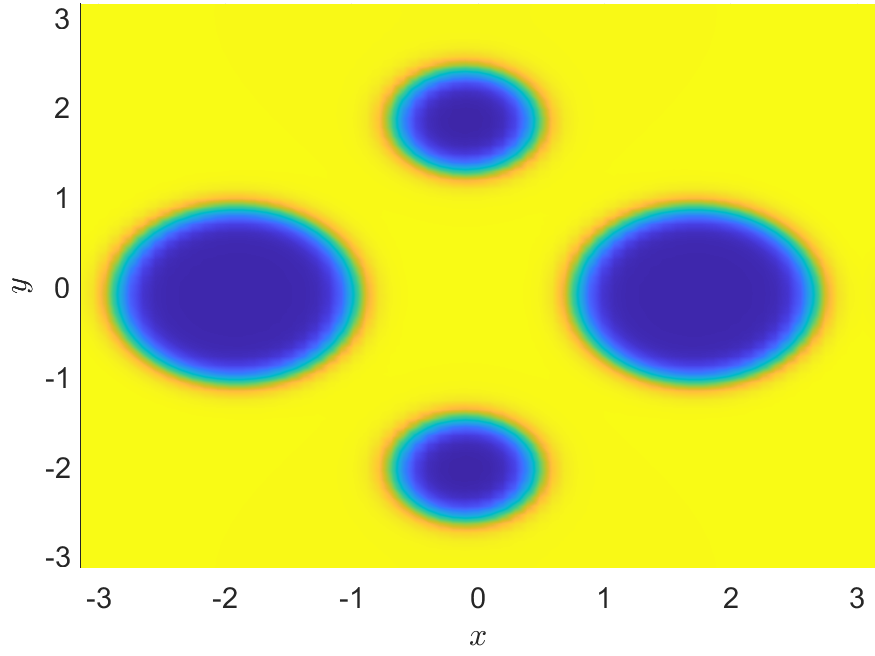}}
	\subfigure{
		\includegraphics[width=1.8in]{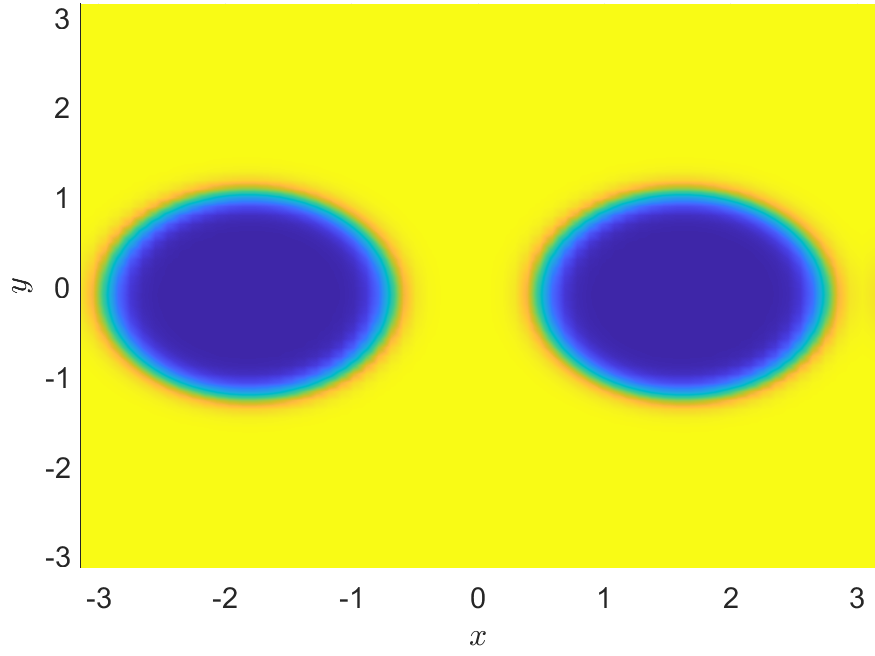}}
	\subfigure{
		\includegraphics[width=1.8in]{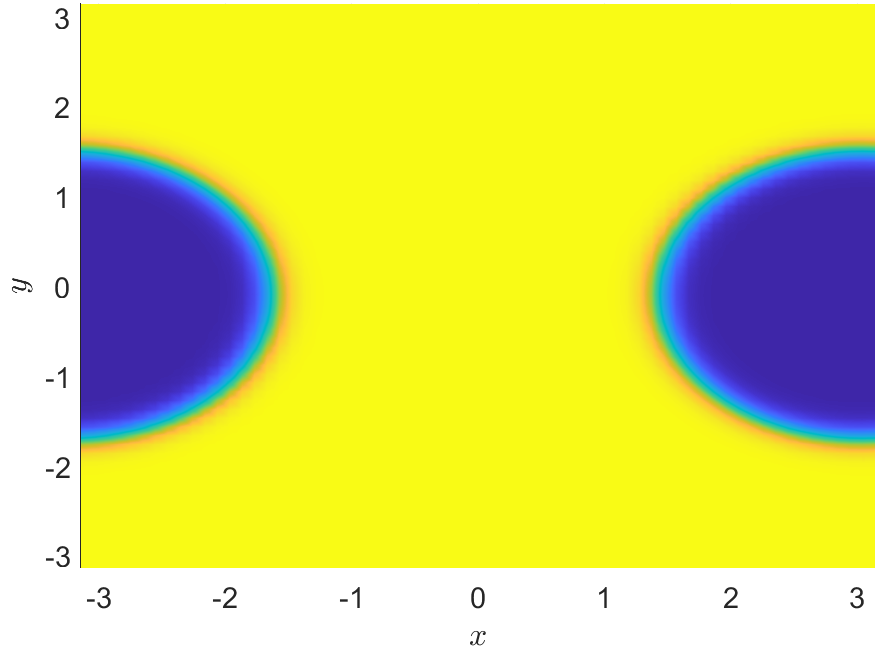}}
	\subfigure{
		\includegraphics[width=1.8in]{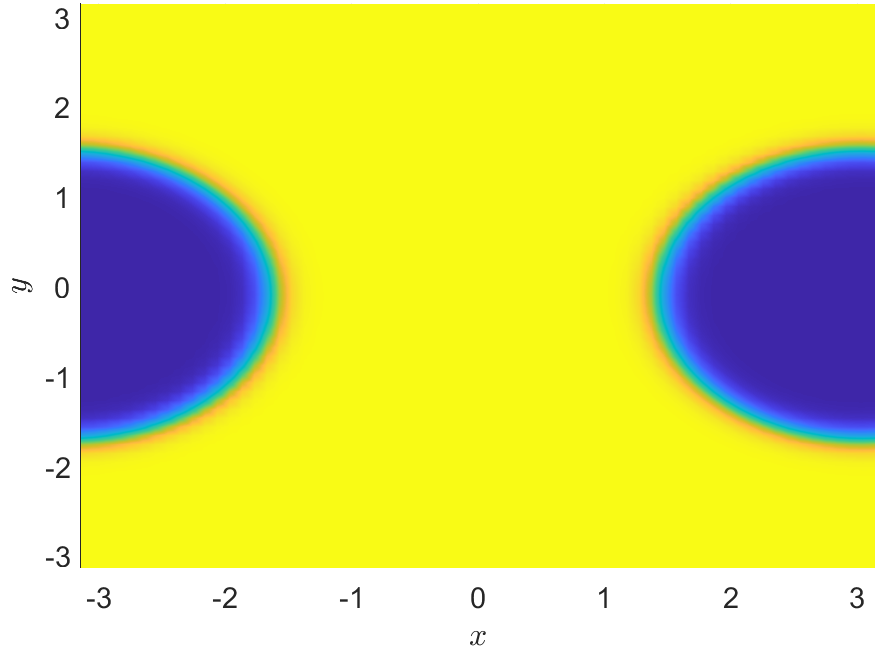}}
	\caption{Solution profiles at $t=0, 2, 20, 200, 450$ and 1000 generated by the R-IERK(3,6;$\tfrac34$) method using the adaptive time-stepping with $\eta=500$ and $\tau_{\max}=0.5$.}
	\label{0fig: IERK decay a52}
\end{figure}

The evolution of micro\/structure for the CH model due to the phase separation at different time are summarized in Figure \ref{0fig: IERK decay a52}, where the phase profiles at $t=0, 2, 20, 200, 450$ and 1000 generated by the R-IERK(3,6;$\tfrac34$) method using the adaptive time-stepping with $\eta=500$ and $\tau_{\max}=0.5$ are depicted. As seen, the micro\/structure is relatively fine and consists of many precipitations at early time and the coarsening, dissolution, merging processes are observed in approaching the steady state.

\section{Concluding remarks}

We propose a class of R-IERK methods in which the associated differentiation matrices and the average dissipation rates are always independent of the time-space discretization meshes. Numerical tests suggest that the R-IERK methods have significant robustnesss in self-adaptive time-stepping procedures as some larger adaptive step-sizes in actual simulations become possible. With the help of an updated time-space error splitting approach with discrete orthogonal convolution kernels and the Gr\"{o}nwall-type lemma for multi-stage methods, the uniform boundedness of stage solutions is theoretically verified without the global Lipschitz continuity assumption of nonlinear bulk so that one can establish the original energy dissipation laws at discrete time stages. Our analysis paves a new way to the internal nonlinear stability of some efficient (not necessarily algebraically stable) Runge-Kutta methods for semilinear parabolic problems. Actually it is closely related to the so-called internal stability analysis \cite{KennedyCarpenter:2003,KennedyCarpenter:2016}, which will be useful to control the stability associated with each stage in addition to each step beyond traditional step-wise stability. 




\end{document}